\documentclass[reqno,a4paper,12pt]{amsart}
\usepackage[margin=3cm,footskip=1cm]{geometry}

\usepackage[latin1]{inputenc}
\usepackage[T1]{fontenc}
\usepackage[english]{babel}
\usepackage{amsmath,amssymb,amscd}
\usepackage{hyperref}
\usepackage{mathtools}
\mathtoolsset{centercolon}

\newcommand\numbersZ{\mathbb{Z}}
\newcommand\numbersR{\mathbb{R}}
\newcommand\numbersC{\mathbb{C}}

\newcommand\lineBundle{\mathcal{L}}
\newcommand\Divisor{\mathfrak{d}}
\newcommand\Deligne{\mathcal{D}}
\newcommand\inverse{^{-1}}
\newcommand\covering{\mathcal{U}}
\newcommand\setC{\check{C}}
\newcommand\numbersRmodZ{\numbersR/\numbersZ}
\DeclareMathOperator{\Map}{Map}
\DeclareMathOperator\Image{Im}
\newcommand\simplex{\Delta}
\DeclarePairedDelimiter{\norm}{\lVert}{\rVert}
\DeclarePairedDelimiter{\abs}{\lvert}{\rvert}
\DeclarePairedDelimiter{\innerp}{\langle}{\rangle}
\newcommand\superprime[1][p]{^{\smash{(#1)}}}

\usepackage[all,cmtip,line,dvips]{xy}

\newtheorem*{abelstheorem}{Abel's theorem}
\newtheorem*{modulitheorem}{Moduli theorem}

\theoremstyle{plain}
\newtheorem{theorem}{Theorem}[section]
\newtheorem{proposition}[theorem]{Proposition}
\newtheorem{lemma}[theorem]{Lemma}
\newtheorem{corollary}[theorem]{Corollary} 
 
\theoremstyle{definition}
\newtheorem{definition}[theorem]{Definition}
\newtheorem{remark}[theorem]{Remark}
\newtheorem{remarks}[theorem]{Remarks}
\newtheorem{example}[theorem]{Example}

\theoremstyle{remark}

\numberwithin{equation}{section}

\def\paritem#1{%
  \smallskip
  \setbox0=\hbox{#1\enspace}
  \par\noindent
  \ifnum\wd0>\parindent\box0
  \else\hbox to\parindent{\box0\hfill}\fi\ignorespaces}

\def\N{\mathbb{N}}
\def\Z{\mathbb{Z}}
\def\R{\mathbb{R}}

\def\D{\mathcal{D}}

\def\H{\mathcal{H}}
\def\I{\mathcal{I}}
\def\L{\mathcal{L}}
\def\M{\mathcal{M}}

\def\U{\mathcal{U}}
\def\V{\mathcal{V}}

\def\epsilon{\varepsilon}

\def\pdd#1#2{\frac{\partial{#1}}{\partial{#2}}}

\DeclareMathOperator\Hom{Hom}

\DeclareMathOperator{\id}{id}
\DeclareMathOperator{\supp}{supp}

\newcommand{\wti}{\widetilde}

\newcommand{\Pic}{\operatorname{Pic}}
\newcommand{\J}{\operatorname{Jac}}

\newcommand{\Tor}{\operatorname{Tor}}

\newcommand{\Pf}{\operatorname{Pf}}
\newcommand{\Deg}{\operatorname{Deg}}
\newcommand{\Star}{\operatorname{Star}}
\newcommand{\Vol}{\operatorname{Vol}}

\newcommand{\veps}{\varepsilon}

\newcommand{\el}{\ell}

\newcommand{\ra}{\rightarrow}
\newcommand{\xra}{\xrightarrow}

\newcommand{\lra}{\longrightarrow}

\newcommand{\x}{\times}

\newcommand{\usetl[1]}{\underset{#1} {\lra}}

\newcommand{\usetlleft[1]}{\underset{#1} {\lla}}

\begin{document}



\title[A generalization of Abel's Theorem ...]
{A generalization of Abel's Theorem \\
and the Abel--Jacobi map${}^1$}
\author[J.~L.~Dupont]{Johan L.~Dupont}
\address{Department of Mathematical Sciences \\
University of Aarhus \\
DK-8000 {\AA}rhus C, Denmark}
\email[J.~L.~Dupont]{dupont@imf.au.dk}
\author[F.~W.~Kamber]{Franz ~W.~Kamber}
\address{Department of Mathematics \\ 
University of Illinois at Urbana--Champaign \\
1409 W. Green Street \\
Urbana, IL 61801, USA}
\email[F.~W.~Kamber]{kamber@math.uiuc.edu}


\thanks{
${}^1$ Work supported in part by the Erwin Schr\"odinger International 
Institute of Mathematical Physics, Vienna, Austria and by the 
Forskningsr{\aa}det for Naturvidenskab og Univers, Denmark
} 


\date{\today}

\begin{abstract}
We generalize Abel's classical theorem on linear equivalence of 
divisors on a Riemann surface. 
For every closed submanifold $M^d \subset X^n$ in a compact oriented 
Riemannian $n$--manifold, or more generally for any $d$--cycle $Z$ 
relative to a triangulation of $X$, we define a (simplicial) 
$(n-d-1)$--gerbe $\Lambda_{Z}$, the Abel gerbe determined by $Z$, 
whose vanishing as a Deligne cohomology class generalizes the notion 
of `linear equivalence to zero'. 
In this setting, Abel's theorem remains valid. 
Moreover, we generalize the classical Inversion Theorem for the Abel--Jacobi 
map, thereby proving that the moduli space of Abel gerbes is isomorphic 
to the harmonic Deligne cohomology; that is, gerbes with harmonic 
curvature. .  
\end{abstract}

\maketitle

\tableofcontents

\setcounter{section}{0}

\section{Introduction}
\label{secone}
In this paper we shall expand on some beautiful ideas of Hitchin \cite{H} and 
Chatterjee \cite{Ch}, generalizing the classical notion of linear equivalence 
of divisors and Abel's theorem about the existence of meromorphic functions 
with prescribed zeroes and poles on a compact Riemann surface (see 
Section \ref{sectwo}). As is well-known, this problem is equivalent to the 
existence of a parallel section, for some complex connection, in the holomorphic 
line bundle of the divisor. 
In general, for a closed oriented Riemannian manifold $X$ of dimension $n$, we
replace the divisor by a cycle $Z$ of arbitrary dimension $d, ~d=0,\dots,n-1$ for 
a smooth triangulation of $X$. 

In Section \ref{secfour}, we associate to $Z$ an abelian \textit{gerbe} 
$\Lambda_{Z}$ which we call the \textit{Abel gerbe} for $Z$, whose class 
$[\Lambda_{Z}]$ in the smooth Deligne cohomology 
$H_{\D}^{n-d} (X, \numbersZ)$ only depends on $Z$. 
Two cycles are then defined to be \textit{linearly equivalent}, if their 
Abel gerbes represent the same class in Deligne cohomology. This definition 
is in agreement with the definition in the classical situation. 
At this level of generality we prove in Section \ref{secfive} Abel's theorem 
\ref{abelone}, characterizing linear equivalence of Abel gerbes in terms of 
period integrals 
(cf. Chatterjee \cite{Ch}, Theorem 6.4.2 for $2$--gerbes associated to
submanifolds of codimension $3$). 

\begin{abelstheorem}
Let $Z=\partial\Gamma$, $\Gamma\in C_{d+1}(K)$. Then $Z$ is linearly
  equivalent to zero, that is $[\Lambda_{Z}] = 0 \in H_{\D}^{n-d}(X, \Z)$, 
  if and only if 
\begin{equation*}   
  \int_\Gamma \theta\in\numbersZ ~,
\end{equation*}
  for all harmonic $\theta \in \Omega^{d+1}(X)$ with integral periods.
\end{abelstheorem}
A similar result was proved by Harvey--Lawson~\cite{HL1} in terms of
`sparks' which however only indirectly provide `gerbes' in the usual 
\v{C}ech--deRham complex for Deligne cohomology 
(see \cite{HL2} for this relation). 

Other well-known results from the theory of Riemann surfaces 
make sense in higher dimensions as well. 
Thus in Section \ref{secsix}, we introduce the \textit{Picard torus} 
of Deligne classes represented by \textit{topologically trivial flat} 
gerbes and the \textit{Jacobi torus} which is the recipient of the 
period map. 
The former is analogous to the Picard variety of holomorphic line 
bundles of degree zero on a Riemann surface, in which case every 
holomorphic line bundle is associated to a divisor. 
The Jacobi torus is analogous to the Jacobi variety of a Riemann 
surface. 
We determine the moduli space $\M_{d} (X)$ of Abel gerbes in 
full generality, as well as the moduli space $\M_d^{\circ} (X)$ 
of topologically trivial Abel gerbes. 

Prior to stating and proving our main theorem \ref{modulitheorem}, 
we illustrate our method by a number of examples 
(Examples \ref{example1} to \ref{example7}). 
Below, we quote Theorem \ref{modulitheorem}.  

\begin{modulitheorem}
  Let $X$ be a compact connected oriented Riemannian manifold $X$ of 
  dimension $n \geq 2$ and let $d = 0,\dots,n-1$. Then  

\paritem{$(1)$} 
  The Picard map $\overline{\alpha} \colon \M_d^{\circ} (X) \to \Pic^{n-d-1}(X)$ 
  is an isomorphism. 
  
\paritem{$(2)$} 
  The Abel--Jacobi map $\overline{J} \colon \M_d^{\circ} (X) \to \J^{d+1}(X)$ 
  is an isomorphism.  
  
\paritem{$(3)$} 
  The mapping $\overline{\Lambda} \colon \M_d (X) \to \H_{\D}^{n-d} (X, \Z)$ 
  is an isomorphism. 
 
\paritem{$(4)$}
  Every equivalence class of $(n-d-1)$--gerbes in the harmonic Deligne 
  cohomology $\H_{\D}^{n-d} (X, \Z)$, given by classes in 
  $H_{\D}^{n-d} (X, \Z)$ whose curvature is harmonic,  
  can be realized by a unique (up to linear equivalence) Abel gerbe. 
  
\end{modulitheorem}
 
In the final Sections \ref{secseven} and \ref{seceight}, we shall investigate 
the Abel gerbe associated to the fundamental cycle of an embedded closed 
submanifold $M \subset X$. In particular, we shall compare the restriction of 
this gerbe to $M$ with the \textit{characteristic gerbe} (\cite{DK2}) called 
the \textit{Euler gerbe}, which respresents the \textit{Cheeger--Chern--Simons class} 
for the normal bundle with the Riemannian connection and is defined in terms of the 
Pfaffian polynomial. We prove in Theorem \ref{euler1} that these two gerbes 
differ by a third canonical gerbe, called the \textit{difference gerbe}. 
This is a topologically trivial gerbe whose curvature is the difference between 
the harmonic form representing the Poincar\'e dual of $[M] \in H_{d}(X)$ and a 
specific choice for the form representing the Thom class of the normal bundle. 

For the construction of these gerbes it is important to use the representation 
of Deligne cohomology and gerbes by \textit{simplicial differential forms} 
as developed in our previous paper \cite{DK2}. 
For completeness, we recall in Section \ref{secthree} the basic definitions 
and properties of these topics. 

We thank James Glazebrook, Jouko Mikkelson and 
J{\o}rgen Tornehave for helpful discussions. 
Both authors thank the Erwin Schr\"odinger International 
Institute for Mathematical Physics (ESI), Vienna, Austria and the 
seond named author thanks the Department for Mathematical Sciences 
(IMF) at Aarhus University for hospitality and support during the 
preparation of this work. 


\section{Abel's Theorem on linear equivalence 
of divisors on a Riemann surface
\label{sectwo}}

  For motivation, let us recall the classical Abel theorem. 
  Let $X$ be a compact Riemann surface and 
  $\Divisor = \sum_{i=1}^k a_ip_i$, $a_i\in\numbersZ$, 
  $p_i\in X$ a divisor. A  first necessary condition for finding a 
  \textit{meromorphic} function with \textit{zeros} and 
  \textit{poles} exactly in $\{p_i\}$ 
  of order $a_i$, is that the degree 
  $\Deg {\Divisor} = \sum_{i} a_i = 0 \in \mathbb{Z}$; 
  that is, there is chain $\Gamma$ with $\partial\Gamma = \Divisor$.

\begin{abelstheorem} 
Suppose that ${\Deg} (\Divisor) = 0$ and $\Divisor = \partial \Gamma$, 
where $\Gamma$ is a (smooth) $1$--chain on $X$. 
Then $\Divisor$ admits a global meromorphic function, that is 
$\Divisor \sim 0$, if and only if 
$$
\int_{\Gamma} \theta \in \mathbb{Z}~, 
$$
for every harmonic $1$--form $\theta \in \mathcal{H}^{1} (X, \mathbb{Z})$ 
with integral periods. 
\end{abelstheorem}

The relationship with smooth 
connections in the holomorphic line bundle $\L (\Divisor)$ for the 
divisor $\Divisor$ is given by the following Lemma. 

\begin{lemma}
  \label{parallelsection}
$\mathcal{L} (\Divisor)$ admits a non--vanishing holomorphic section; 
that is, $\Divisor \sim 0$, if and only if $\mathcal{L} (\Divisor)$ 
admits a non--vanishing $C^{\infty}$--section, which is parallel with 
respect to a suitable complex 
connection in $\mathcal{L} (\Divisor)$. 
\end{lemma}

\begin{proof}
Recall that $\mathcal{L} (\Divisor)$ is given by first choosing 
an open covering $\covering=\{U_i\}_{i\in I}$ and local solutions 
$f_i$ on $U_i$. Then
\begin{equation*}
  g_{ij}=f_i/f_j\colon U_i\cap U_j \to \numbersC^*=\numbersC\setminus\{0\}
\end{equation*}
is a \textit{\v{C}ech cocycle} defining a \textit{holomorphic} line bundle
$\lineBundle(\Divisor)$. If $h_i \colon U_i \to \mathbb{C}^{\ast}$ defines 
a \emph{holomorphic} section of $\lineBundle(\Divisor)$
i.e.\ if $g_{ij}h_j=h_i$ on $U_i\cap U_j$ $\forall
i,j\in I$, then
\begin{equation*}
  f_i/h_i = f_i / (g_{ij}h_j) = f_j/h_j
\end{equation*}
defines a \textit{global meromorphic solution}. For finding $\{h_i\}$ we
define a smooth connection in $\lineBundle(\Divisor)$, i.e.\ a
family $\omega_i\in\Omega^1(U_i)$, such that \ $g_{ij}\inverse d g_{ij}
=\omega_i-\omega_j$ on $U_i\cap U_j$ and we can arrange that
$\omega_i=f_i\inverse df_i$ away from small neighborhoods of $\{p_i\}$.

Now suppose $\lineBundle(\Divisor)$ has a non--vanishing \textit{parallel}
    $C^\infty$--section $\{k_i\}$, i.e.\ a section satisfying
  \begin{equation*}
    k_i\inverse dk_i=\omega_i\quad \text{in $U_i$}\ .
  \end{equation*}
  Then away from $\{p_i\}$ we have
  \begin{equation*}
    d\log k_i = \omega_i = d\log f_i\ .
  \end{equation*}
  Hence $\log k_i$ is holomorphic away from $p_i$. But $\log k_i$ is
  smooth all through $U_i$ so the singularity of $\log k_i$ is
  removable. Hence we can redefine $k_i$ throughout $U_i$ to give a
  holomorphic section.
\end{proof}

Our goal is to generalize these classical results to submanifolds 
$M^{d} \subset X^{n}$ of compact oriented Riemannian manifolds $X^{n}$, 
and more generally to cycles $Z \subset X$, by using the notion of the
Abel gerbe. 


\medbreak
\section{Review of `gerbes with connection' and simplicial gerbes}
\label{secthree}

\subsection{Gerbes with connections}
\label{subsecthreeone}

Let $X$ be a smooth manifold and $\covering =\{U_i\}_{i\in I}$ an open
covering. We assume throughout that the covering $\covering$ is good, 
i.e.\ all
\begin{equation*}
  U_{i_0 \dots i_p}=U_{i_0} \cap \dots \cap U_{i_p}
\end{equation*}
are contractible.
We identify the circle group $U (1)$ with $\numbersRmodZ$ via the 
exponential map; that is
\begin{align*}
  U(1) = \text{circle group} &\cong \numbersRmodZ\ . 
  \\
  \exp(2\pi it)   &\leftrightarrow t 
\end{align*}
A \textit{Hermitian} $\el$--\textit{gerbe} 
is given by a cocycle in the \v{C}ech complex   
\begin{equation*}
  \theta\in \setC^{\el} (\covering,\underline{\numbersRmodZ})\ ;
\end{equation*}
that is, $\theta_{i_0\dots i_{\el} }\colon U_{i_0\dots i_{\el}} 
\to \numbersRmodZ$ satisfying 
\begin{equation*}
  0\equiv \check{\delta}\theta_{i_0\dots i_{\el} }=\sum_\nu(-1)^\nu
  \theta_{i_0\dots \check{i}_\nu\dots i_{\el+1} } \ .
\end{equation*}
For $\el=1$, $\theta$ defines a line bundle.

We consider the modified \v{C}ech-deRham bi--complex:
\begin{equation*}
  \xymatrix{
    & \vdots & \vdots & \vdots &
    \\
    \setC^2(\covering,\numbersRmodZ) \ar[r] 
    &
    \setC^2(\covering,\underline{\numbersRmodZ})\ar[r]^-d \ar[u] 
    &
    \setC^2(\covering,\underline{\Omega}^1) \ar[r]^-d \ar[u] \ar@{.}[ul]
    &
    \setC^2(\covering,\underline{\Omega}^2)\ar[r]^-d \ar[u] \ar@{.}[ul]
    &
    \cdots
    \\
    \setC^1(\covering,\numbersRmodZ) \ar[r] \ar[u]_\delta
    &
    \setC^1(\covering,\underline{\numbersRmodZ})\ar[r]^-d \ar[u] \ar[u]_\delta
    &
    \setC^1(\covering,\underline{\Omega}^1) \ar[r]^-d \ar[u] \ar@{.}[ul]\ar[u]_\delta
    &
    \setC^1(\covering,\underline{\Omega}^2)\ar[r]^-d \ar[u] \ar@{.}[ul]\ar[u]_\delta
    &
    \cdots
    \\
    \setC^0(\covering,\numbersRmodZ) \ar[r] \ar[u]_\delta
    &
    \setC^0(\covering,\underline{\numbersRmodZ})\ar[r]^-d \ar[u] \ar[u]_\delta
    &
    \setC^0 (\covering,\underline{\Omega}^1) \ar[r]^-d \ar[u] \ar@{.}[ul]\ar[u]_\delta
    &
    \setC^0(\covering,\underline{\Omega}^2)\ar[r]^-d \ar[u] \ar@{.}[ul]\ar[u]_\delta
    &
    \cdots
    \\
    &
    \Map (X,\underline{\numbersRmodZ}) \ar[r]^-d \ar[u]_{\epsilon^*}
    &
    \Omega^1(X) \ar[r]^-d \ar[u]_{\epsilon^*}
    &
    \Omega^2(X) \ar[r]^-d \ar[u]_{\epsilon^*}
    &
    \cdots
  }
\end{equation*}
where the dotted lines indicate the total complex with differential 
$D=\check\delta+(-1)^p d$ on $\check{C}^p (\covering,\Omega^*)$. 

A \textit{connection} in an $\el$--gerbe $\theta$ is a 
sequence
$\omega=(\omega^0,\dots,\omega^{\el})$ in the \v{C}ech--deRham
bi\nobreakdash-complex
\begin{equation*}
  \omega^\nu \in \setC^\nu(\covering,\Omega^{\el-\nu}), \quad \nu =0,\dots,\el\ , 
\end{equation*}
satisfying
\begin{equation*}
  \omega^{\el} \equiv -\theta \bmod{\numbersZ}\ , \quad
  \check{\delta} \omega^{\nu-1}+(-1)^\nu d\omega^\nu=0\ , \quad
  \nu=1,\dots,\el\ .
\end{equation*}
In particular, we have $\check\delta (d \omega^0) = 0$, so that $d \omega^0$ 
is given by a global form $F_{\omega}$, the \textit{curvature} of 
$(\theta, \omega)$; that is, we set 
\begin{align*}
  d\omega^0 & \in \Image\bigl\{\epsilon^*\colon
  \Omega^{\el+1}(X) \to \setC^0 (\covering,\Omega^{\el+1})\bigr\}\ , 
  \\
  \text{where} \quad \epsilon &\colon \sqcup_i U_i\to X \quad 
  \text{is the natural map and}
  \\
  F_\omega &:= (\epsilon^*)\inverse(d\omega^0) \in \Omega^{\ell +1} (X)\ . 
\end{align*}

\begin{definition}
\paritem{($1$)}  
  Two gerbes with connection are \textit{equivalent}, 
  $(\theta_1,\omega_1)\sim (\theta_2,\omega_2)$,  
  if $\omega_1-\omega_2$ is a coboundary in
  \begin{equation*}
    \bigl(\setC^*(\covering,\Omega^*)\big/\setC^*(\covering,\numbersZ),
    D \bigr)\ . 
  \end{equation*}

\paritem{($2$)}
$H_{\Deligne}^{\el+1}(X,\numbersZ)$, the set of equivalence classes 
$[\theta, \omega]$ of $\el$--gerbes with connection, is the smooth 
Deligne cohomology of $X$. 

\paritem{($3$)}
$H^{\el} (X,\numbersRmodZ)$ is the set of equivalence classes of
$\el$-\nobreakdash-gerbes with \textit{flat} connection; that is 
$F_\omega=0$. Hence we have the exact sequence
  \begin{equation}
  \xymatrix@R=0mm{
        0 \ar[r] 
        &
        H^{\el} (X,\numbersRmodZ) \ar[r]
        &
        H_{\Deligne}^{\el+1}(X,\numbersZ)\ar[r]^{d_{\ast}} 
        &
        \Omega_{\mathrm{cl}}^{\el+1}(X,\numbersZ) \ar[r]
        &
        0\ ,
        \\
        & &
        [\theta,\omega] \ar@{|->}[r]
        &
        F_\omega
      }
  \label{DeligneES}
  \end{equation}
    where $\Omega_{\mathrm{cl}}^{\el+1}(X,\numbersZ)$ denotes the closed
    $(\el+1)$--forms with \textit{integral periods}. 

\end{definition}


Let us introduce the notation
\begin{align}
  \label{global}
  H^{\ell+1}_\D(X)=\Omega^\ell (X)/ d\Omega^{\ell-1}(X)\ . 
\end{align}
The elements $[\omega] \in H^{\ell+1}_\D (X)$ can be interpreted as equivalence 
classes of connections on the trivial $\ell$--gerbe $\theta = 0$ by setting 
\begin{equation}
  \label{globalDeligne}
  \omega^0=\veps^* \omega\ , \qquad 
  F_{\omega}= d\omega\ , \qquad 
  \check{\delta} \omega^0=0\ , \qquad 
\omega^1=\ldots =\omega^{\ell} = 0\ . 
\end{equation} 
Thus $\iota (\omega ) = (0; \veps^{\ast}\omega, 0, \dots, 0)$ 
induces a well--defined mapping 
\begin{equation*}
\iota_{\ast} \colon H^{\ell+1}_\D(X) \ra H_\D^{\ell+1}(X,\Z)\ , 
\end{equation*}
since $\iota (d\alpha ) = D (0; \veps^{\ast}\alpha, 0, \dots, 0)$. 
Clearly the connection is flat if and only if $F_\omega = d\omega=0$, 
that is $[\omega] \in H^{\ell} (X,\R)$. 

\vfill\newpage
We then have the following commutative diagram with exact rows and columns: 
\begin{equation}
   \label{DeligneDiagram}
\vcenter{
  \xymatrix{
& 0 \ar[d] & 0 \ar[d]  \\
0 \ar[r] & j_{\ast} H^\ell (X,\Z) \ar[d] \ar[r]^{\cong} & 
\Omega_{\rm cl}^{\ell} (X ,\Z) / d \Omega^{\ell -1}(X)
\ar[d] \ar[r] & 0 \ar[d] \\ 
0 \ar[r] & H^\ell (X,\R) \ar[d]^{\rho_{\ast}} \ar[r] & H_\D^{\ell+1} (X) 
\ar[d]^{\iota_{\ast}} \ar[r] &
\Omega^\ell (X) / \Omega_{\rm cl}^{\ell} (X) \ar[r] \ar[d]^d  &0 \\  
0 \ar[r] & H^\ell(X,\R / \Z) \ar[d]^{\beta_{\ast}} \ar[r] & H_\D^{\ell+1}(X,\Z) \ar[d] 
\ar@{->>}[dl]_{c} \ar[r]^{d_*} &\Omega_{\rm cl}^{\ell+1}(X,\Z) \ar[d]
\ar[r] &0   
\\
& H^{\ell +1}(X,\Z) \ar[d]^{j_{\ast}}  
& H^{\ell}(X,\underline{\R / \Z}) \ar[d] |!{[r];[dl]}\hole
\ar[l]_{{\check\delta}_{\ast}}^{\cong} 
& j_{\ast} H^{\ell +1} (X,\Z) \ar@{>->}[dll] \ar[d] 
& 0  \ar[l]   
\\
& H^{\ell +1}(X,\R)  & 0 & 0 
}
}
\end{equation}

\begin{remarks}
  \label{delignerem} 

The diagram \eqref{DeligneDiagram} incorporates many 
properties of our construction: 

	
	\paritem{(1)} 
	The second exact row follows from the definition \eqref{global}. 
	
	\paritem{(2)} 
	The third exact row is \eqref{DeligneES}, with $d_{\ast}$ being 
	the curvature. 

	\paritem{(3)}
  The exact column on the left is the Bockstein sequence for 
  \begin{equation*}
  0 \ra \numbersZ \xra{j} \numbersR \xra{\rho} \numbersRmodZ \ra 0 . 
  \end{equation*}
	Note that the image of $\beta_{\ast}$ is the (finite) torsion subgroup of 
	$H^{\ell +1}(X,\Z)$; that is, $\beta_{\ast}$ induces an isomorphism 
  \begin{equation*}
  H^{\ell}(X, \R / \Z) / 
	\overline{\rho}_{\ast} \left( H^\ell (X,\R) / j_{\ast} H^\ell (X,\Z) \right) 
	\cong \Tor_{\Z} (H^{\ell +1} (X,\Z), \numbersRmodZ) \subseteq H^{\ell +1} (X,\Z)\ . 
	\end{equation*}
	
	\paritem{(4)}
	The map $c$ is the characteristic class 
	$\check\delta_{\ast} [\theta] = - [\check\delta \omega^{l}]$ of the gerbe 
	$[\theta,\omega]$ (the Douady--Dixmier invariant); it is equivalent to the 
	last map $[\theta,\omega] \mapsto [\theta]$ in the middle exact column, 
	which simply forgets the connection. These maps are surjective, since every 
	(naked) gerbe $[\theta] \in H^{\ell}(X,\underline{\R / \Z})$ admits a
	connection. 
	
	\paritem{(5)}\label{item325}
 	The image of $\iota_{\ast}$~; that is, the equivalence classes of 
 	trivial gerbes with connection, is given exactly by the 
 	the kernel of the characteristic class $c$, so we may call 
 	these gerbes \textit{topologically trivial}. 
	
	\paritem{(6)}
	It follows that the Deligne cohomology is given by an exact sequence 
	(i.e. the middle exact column in \eqref{DeligneDiagram})
	\begin{equation}
    0 \ra \Omega^{\ell} (X)  / \Omega_{\rm cl}^{\ell} (X ,\mathbb{Z})
    \xra[\subseteq]{\overline{\iota}_{\ast}} 
    H_\D^{\ell+1}(X,\Z) \xra{c} H^{\ell +1}(X,\Z) \ra 0 ~. 
  \label{DeligneES2}
  \end{equation}
  
  \paritem{(7)}
	The commutativity of the diagram involving the slanted arrows expresses the fact 
	that the characteristic class of a gerbe determines the deRham class of the 
	curvature in $H^{\ell +1}(X,\R)$; that is, 
	$j_{\ast} c ([\theta,\omega]) = [F_{\omega}]$.

\end{remarks}


\subsection{Simplicial forms and gerbes}
\label{subsecthreetwo}

In this section, we recall the reformulation of Deligne
cohomology in terms of simplicial deRham theory \cite{DK2}. 
For simplicial deRham theory we refer to \cite{D1}, \cite{D2}. 


Consider the standard simplex $\simplex^p\subseteq \numbersR^{p+1}$
\begin{gather*}
  \simplex^p = \{(t_0,\dots,t_p)\mid{\textstyle \sum_i } t_i=1, \ t_i
  \geq 0\}
  \\
  \text{with face maps} \quad
  \epsilon^i\colon \simplex^{p-1} \to \simplex^p\ , i=0,\dots,p\ , 
  \text{ given by}
  \\
  \epsilon^i(t_0,\dots,t_{p-1}) = (t_0,\dots,0,\dots,t_{p-1}), ~
  (t_0,\dots,t_{p-1}) \in \simplex^{p-1}\ .
\end{gather*}
The open covering $\covering =\{U_i\}_{i\in I}$ of $X$ determines a 
\textit{simplicial manifold} $N\covering$
\begin{equation*}
  N\covering(p) = \bigsqcup_{(j_0,\dots, j_p)} U_{j_0\dots j_p}, \quad
  p=0,1,\dots 
\end{equation*}
\begin{align*}
  \text{with} ~\textit{face operators} \quad 
  \epsilon_j\colon N\covering(p) &\to N\covering(p-1)\ , 
  \quad i=0,\dots,p\ , \quad \text{given by}
  \\
  U_{j_0\dots j_p} &\hookrightarrow U_{j_0\dots\check{j}_i\dots j_p}\ . 
\end{align*}
The \textit{fat realisation} is 
\begin{equation*}
  \norm{N\covering} = \bigsqcup_p \simplex^p\times N\covering(p)/\sim ,
\end{equation*}
with identifications 
$(t,\epsilon_i x)\sim(\epsilon^i t,x)$, $t\in\simplex ^{p-1}$, 
$x\in N\covering (p)$.

\begin{definition}
  A \textit{simplicial} $k$--form $\omega$ on $N\covering$ is a
  sequence $\omega\superprime \in
  \Omega^k\bigl(\simplex^p\times N\covering(p)\bigr)$ satisfying
  \begin{equation*}
    (\epsilon^i\times\id)^*\omega\superprime = (\mathord{\id}\times
    \epsilon_i)^*\omega\superprime[p-1],\quad
    i=0,\dots,p,\ \forall p ~, 
  \end{equation*}
\end{definition}
\noindent
and we denote by $\Omega^k(\norm{N\covering})$ the set of simplicial
$k$--forms. 

\medbreak
Recalling that the open covering $\covering$ is assumed to be good, 
we have the following results. 
\begin{theorem}[deRham] \cite{D1}, \cite{D2} 
  \label{derham}
  There are quasi--isomorphisms 
  (inducing isomorphisms in cohomology)
  \begin{equation*}
    \I_\simplex \colon \Omega^*(\norm{N\covering}) \to
    \setC(\covering, \Omega^*)\ , 
  \end{equation*}
  given by
  \begin{equation*}
    \I_\simplex(\omega)=(\omega^\nu)\ , \qquad \omega^\nu =
    \int_{\simplex^\nu} \omega\superprime[\nu]\ ;
  \end{equation*}
  and
  \begin{equation*}
  \epsilon^* \colon \Omega^*(X)\to \Omega^*(\norm{N\covering})\ , 
  \end{equation*}
  induced by the natural map
  \begin{equation*}
    \epsilon \colon \simplex^p\times N\covering(p)\to N\covering(p) \to X\ . 
  \end{equation*}
\end{theorem}
We also need the following 
\begin{definition}
  $\omega\in\Omega^k(\norm{N\covering})$ is \textit{integral} if
  \begin{enumerate}
  \item $\omega\superprime = \sum \alpha_{i_0, \dots , i_k}(t) ~
  dt_{i_1}\wedge \dots \wedge dt_{i_k}$\ , 
  \item $\I_\simplex(\omega)\in \setC^*(\covering,\numbersZ)\subseteq
    \setC^*(\covering, \Omega^0)$\ . 
  \end{enumerate}
We denote by 
$\Omega^*_{\numbersZ}(\norm{N\covering})\subseteq\Omega^*(\norm{N\covering})$
the subcomplex of integral forms. 
\end{definition}

\begin{remark}
  Note that we now have that
  $\I_\simplex\colon \Omega^*_{\numbersZ}(\norm{N\covering}) \to 
  \setC^*(\covering, \numbersZ)$ is also a 
  quasi--isomorphism.
\end{remark}

\begin{theorem} \cite{DK2} 
  \label{simpldeligne1}
  Every $\el$--gerbe with connection $(\theta, \omega)$ 
  is up to equivalence determined by a simplicial form
  $\Lambda\in\Omega^{\el}(\norm{N\covering})$ satisfying
  \begin{equation}
    \label{simpldeligne2}
    d\Lambda=\epsilon^*\alpha-\beta\ , \quad \text{with }
    \alpha\in\Omega^{\el+1}(X)\ , 
    ~\beta\in\Omega^{\el+1}_{\numbersZ}(\norm{N\covering})\ . 
  \end{equation}
  In fact
  \begin{equation*}
    \omega^\nu =\int_{\simplex^{\nu}} \Lambda^{\nu}\ , 
    \quad \nu=0,\dots,\el\ ,\quad -\theta = \omega^{\el}\ ,
  \end{equation*}
  and $\alpha$ is the curvature.
\end{theorem}

Equivalently, we have 
\begin{theorem}\label{simpldeligne3}
  Every element in $H_{\Deligne}^{\el+1}(X,\numbersZ)$ is 
  represented by a unique class $[\Lambda]$ in
  \begin{equation*}
    \Omega^{\el} (\norm{N\covering})\big/\bigl(\Omega^{\el}_{\numbersZ}
    (\norm{N\covering})+d\Omega^{\el-1}(\norm{N\covering})\bigr)\ , 
  \end{equation*}
  satisfying \eqref{simpldeligne2} above.
\end{theorem}

\begin{proof}
  Let $H_{\Deligne}^{\el+1}(X,\numbersZ)'$ be the subgroup of such classes
  $[\Lambda]$ satisfying \eqref{simpldeligne2}. Then there is a diagram 
  with exact rows:
  \begin{equation*}
    \xymatrix@C=6mm{
      0 \ar[r] 
      &
      H^{\el}\bigl(\Omega^*(\norm{N\covering})\big /
      \Omega^*_{\numbersZ}(\norm{N\covering})\bigr)  
      \ar[r] \ar[d]^{\I_\simplex}_{\cong}
      &
      H_{\Deligne}^{\el+1}(X,\numbersZ)' \ar[r] \ar[d]^{\I_\simplex}
      &
      \Omega^{\el+1}_{\mathrm{cl}}(X,\numbersZ)\ar[r] \ar@{-}[d]^{\id}
      & 
      0
      \\
      0 \ar[r]
      & 
      H^{\el} (X,\numbersRmodZ) \ar[r]
      &
      H_{\Deligne}^{\el+1}(X,\numbersZ) \ar[r]
      &
      \Omega^{\el+1}_{\mathrm{cl}}(X,\numbersZ) \ar[r]
      &
      0
    }
  \end{equation*}
  The vertical map on the left is an isomorphism by deRham's Theorem. 
  Hence Theorem~\ref{simpldeligne3} follows from the $5$\nobreakdash-lemma.
\end{proof}


\medbreak
\section{Abel gerbes associated to cycles and submanifolds
\label{secfour}}

Classically on a Riemann surface two divisors $\Divisor_1,\Divisor_2$ are called
\textit{linearly equivalent} if $\Divisor_1-\Divisor_2$ is the divisor of a 
meromorphic function. We have seen that this is equivalent to finding a parallel
section for a suitable connection in the line bundle
$\lineBundle(\Divisor_1- \Divisor_2)$.
Using gerbes we can generalise this to higher dimensions as follows. 

Let $X = X^n$ be a \textit{compact connected oriented} manifold, $\partial
X=\emptyset$, with Riemannian metric. Choose a \textit{smooth 
triangulation}, i.e.\ a homeomorphism to a finite simplicial complex
$X \approx \abs{K}$, such that the homeomorphism is a diffeomorphism on
each simplex.
For a \textit{cycle} 
\begin{equation*}
  Z \in C_d (K)\ , 
\end{equation*}
let $\abs{Z} \subseteq \abs{K}$ be the subcomplex consisting of 
all simplices of $Z$ and their faces. 
For the good covering $\covering$ of $X$ given by the 
\textit{open stars} of the vertices of $K$, we set 
$\covering_{X - \abs{Z}} = \{ U_j \mid U_j \cap \abs{Z} = \emptyset \}$ 
and $\covering_Z = \covering - \covering_{X - \abs{Z}} = 
\{U_i \mid U_i \cap \abs{Z} \neq \emptyset\}$. 
Then $\covering_{X - \abs{Z}}$ is a covering of $X - \abs{Z}$ and 
$\covering_Z$ is a covering of a regular neighborhood of $Z \subset X$. 
Let
\begin{equation*}
  \eta_Z\in \H^{n-d} (X, \Z) \subset \Omega^{n-d}(X) \quad ,\quad 
  \beta_Z\in\Omega_{\Z}^{n-d}(\norm{N\covering}) 
\end{equation*}
both represent the \textit{Poincar\'e dual} of $[Z]\in H_d(X)$: 
$\eta_Z$ is a \textit{harmonic} form with integral periods, 
and $\beta_Z$ is an \textit{integral} form with 
$\supp \beta_{Z} \subseteq \norm{N \covering} - \norm{N \covering_{X - \abs{Z}}}$. 
Here $\mathcal{H}^{\ell} (X, \mathbb{Z} ) \subset \mathcal{H}^{\ell} (X)$ 
denote the harmonic forms, respectively the integral lattice of harmonic forms. 

\medbreak
Following Hitchin~\cite{H}, we can solve the distributional Poisson equation 
in $\Omega^{d}(X)^{\prime}$ ~\cite{DRh}~: 
\begin{equation}
\label{poisson1}
\Delta H_{Z} = \eta_{Z} - \delta_{Z}\ , 
\end{equation}
where $\Delta$ is the Laplace operator, $\eta_{Z}$ the harmonic form dual to $[Z]$ 
and $\delta_{Z}$ the Dirac measure associated to $Z$; that is, 
\begin{equation*}
\eta_{Z} (\psi) = \int_{X} \eta_{Z} \wedge \psi   \quad , \quad 
\delta_{Z}(\psi) = \int_{Z} \psi, ~\psi \in \Omega^{d}(X)\ . 
\end{equation*}
$H_{Z}$ is uniquely defined up to a global harmonic $(n-d)$--form, and is 
smooth outside $\abs{Z}$. Since $\eta_{Z}$ and $\delta_{Z}$ represent the 
same cohomology class, we get from the deRham--Hodge decomposition 
\begin{equation*}
\Delta H_{Z} = d \ast d \ast H_{Z}\ , 
\end{equation*}
where $\ast$ is the Hodge $\ast$--operator. Setting 
$F_{Z} = \ast~d \ast H_{Z}$, it follows that 
$F_{Z}$ is uniquely defined by $Z$ and we have 
\begin{equation}
  \label{poisson2}
\Delta H_{Z} = d F_{Z} = \eta_{Z} - \delta_{Z}\ . 
\end{equation} 
In particular, $F = F_{Z} \vert_{X - \abs{Z}} = \ast~d \ast H_{Z} 
\vert_{X - \abs{Z}}$  is smooth 
and satisfies 
\begin{equation}
  \label{poisson3}
dF = \eta_{Z}\vert_{X - \abs{Z}}~ , ~d^{\ast} F = 0\ . 
\end{equation}

\begin{theorem}
\label{abelgerbeone}
  There is a canonical Deligne class $[\Lambda_Z] \in
  H_{\Deligne}^{n-d}(X,\numbersZ)$, such that 
  $\Lambda_Z \in \Omega^{n-d-1}(\norm{N\covering})$ satisfies: 
  
  \paritem{$(1)$}
  $d\Lambda_Z=\epsilon^*\eta_Z -\beta_Z$. 
  Thus the curvature of $[\Lambda_{Z}]$ is the integral harmonic form 
  $\eta_{Z}\in \H^{n-d}(X, \Z)$ 
  and the characteristic class of $[\Lambda_{Z}]$ is the Poincar\'e dual 
  $[\beta_{Z}] \in H^{n-d} (X, \Z)$ of $[Z] \in H_d (X)$. 

  \paritem{$(2)$}
  $[\Lambda_{Z}]$ is additive; that is, we have 
  $[\Lambda_{Z_1 + Z_2}] = [\Lambda_{Z_1}] + [\Lambda_{Z_2}]$, for 
  $Z_1, Z_2 \in Z_d (K)$. 
  
  \paritem{$(3)$} 
  $F = F_Z \vert_{X - \abs{Z}} \in \Omega^{n-d-1}(X-\abs{Z})$ is smooth and 
  satisfies $\Delta F = 0$; that is, $F$ is harmonic on ${X - \abs{Z}}$. 
  
  \paritem{$(4)$} 
  $\Lambda_{Z} \vert_{W} = \epsilon^*F$, 
  where $W = \norm{N \covering_{X - \abs{Z}}}$. 
\end{theorem}

\begin{proof}
Let 
\begin{equation*}
  \begin{aligned}
  K_{0} &= \{ a_0, \dots, a_m, \dots, a_{N} \} \ ,   \\
  Z_{0} &= \{ a_0, \dots, a_m \}
  \end{aligned}
\end{equation*}
be the vertices of $K$ and the subcomplex $\abs{Z}$ respectively. 
Then the coverings of $\abs{K}$, $\abs{Z}$ and ${X - \abs{Z}}$ 
respectively are given by $\covering = \{ U_i \mid i=0,\ldots, N \}$,  
$\covering_{Z} = \{ U_i \mid i=0,\ldots, m \}$ and  
$\covering_{X - \abs{Z}} = \{ U_j \mid j=m+1 ,\ldots, N \}$, 
where $U_i = \Star(a_i)$. 
Let $V = \bigcup_{i=0}^{m} U_i$, which is a regular neighborhood 
of $\abs{Z}$. Then by Lefschetz and Poincar\'e duality we have 
a commutative diagram 
\begin{equation*}
\xymatrix{
H^{n-d} (\overline{V}, \partial\overline{V} ) \ar[r]^{\cong} \ar[d]^{\cong} & 
H^{n-d} (X, X - \abs{Z}) \ar[r]^{\cong} \ar[d]  & H_d(\abs{Z}) \ar[d]  \\
H^{n-d} (X, X - V) \ar[r] & H^{n-d}(X) \ar[r]^{\cong} & H_d (X) 
} 
\end{equation*}
It follows, as claimed above, that the Poincar\'e dual of $[Z]\in H_d (X)$ 
is represented in $H^{n-d} (X) \cong H^{n-d} (\norm{N \covering})$ by an 
integral simplicial form $\beta_{Z} \in \Omega_{\Z}^{n-d} (\norm{N\covering})$ with 
$\supp \beta_{Z} \subseteq \norm{N \covering} - \norm{N \covering_{X - \abs{Z}}}$. 

For the construction of $\Lambda_{Z}$, we first define the following
simplicial forms 
\begin{equation*}
\eta_{0}, ~\eta_1, ~\eta_2 \in \Omega^{n-d} (\norm{N \covering})
\end{equation*} 
and $F_1 \in \Omega^{n-d-1} )(\norm{N \covering})$. They are given on 
$\simplex^{p} \x U_{i_0} \cap \dots \cap U_{i_p}$ respectively 
by the forms: 
\begin{equation*}
\begin{aligned}
(\eta_{0})_{i_0\dots i_p} &= \sum_{i_s \leq m} t_{i_{s}} ~\eta_{Z} 
\quad , \quad  
(\eta_{1})_{i_0\dots i_p} = \sum_{i_s > m} t_{i_{s}} ~\eta_{Z}\ ,   \\
(\eta_{2})_{i_0\dots i_p} &= \sum_{i_s > m} d t_{i_{s}} \wedge F = - 
\sum_{i_s \leq m} d t_{i_{s}} \wedge F \ ,  \\
(F_{1})_{i_0\dots i_p} &= \sum_{i_s > m} t_{i_{s}} ~F\ . 
\end{aligned}
\end{equation*}
Notice that $\eta_1, \eta_2$ and $F_1$ vanish on
$\simplex^p \x U_{i_0\dots i_p} \cap \abs{Z}$, since 
$U_{i_0\dots i_p} \cap \abs{Z} \neq \emptyset$ only 
if all $i_s \leq m$. 
>From these formulas, we clearly have 
\begin{equation*}
\begin{aligned}
d F_{1} &= \eta_1 + \eta_2  \quad , \quad 
\eta_{Z} = \eta_{0} + \eta_1\ ,  \\ 
\eta_{Z} &= (\eta_{0} - \eta_2) + d F_1\ . 
\end{aligned}
\end{equation*}
The third equation implies that $d (\eta_{0} - \eta_2) = 0$. 
Furthermore, by construction 
\begin{equation*}
\supp \eta_{0} ~, ~\supp \eta_2 \subseteq 
\norm{N \covering} - \norm{N \covering_{X - \abs{Z}}}\ . 
\end{equation*}
It follows that both $\beta_{Z}$ and $\eta_0 - \eta_2$ lie in 
$\Omega^{n-d} (\norm{N \covering})$, both have support in 
$\norm{N \covering} - \norm{N \covering_{X - \abs{Z}}}$ 
and both represent the Lefschetz dual of $[Z]\in H_d (\abs{Z})$ in 
$H^{n-d} (\norm{N \covering}, \norm{N \covering_{X - \abs{Z}}} ) 
\cong H^{n-d}(X, X-V)$. 
Hence there is a simplicial form 
$\gamma \in \Omega^{n-d-1} (\norm{N \covering})$, also with 
$\supp \gamma \subseteq \norm{N \covering} - 
\norm{N \covering_{X - \abs{Z}}}$, such that 
$\eta_{0} - \eta_2 = \beta_{Z}+ d \gamma$. 
Now we define $\Lambda_{Z} = \gamma + F_1$, so that we have 
\begin{equation}
  \label{lambdaequ}
\begin{aligned}
\Lambda_{Z} &= \gamma + F_1 \in \Omega^{n-d-1} (\norm{N \covering})\ ,  \\
d \Lambda_{Z} &= d \gamma + d F_1 = \eta_{Z} - \beta_{Z}\ . 
\end{aligned}
\end{equation}
We must now show that the class of $\Lambda_{Z}$ in Deligne 
cohomology $H_{\D}^{n-d}(X, \Z)$ depends only on $Z$. 
Recalling that $F$ 
is uniquely defined by the Poisson equation \eqref{poisson2}, 
let $\beta^{\prime}$ be another integral form 
representing the Poincar\'e dual of $[Z]$ and suppose 
$\gamma^{\prime}$ satisfies the same properties as 
$\gamma$ relative to $\beta^{\prime}$; in particular 
$\eta_{0} - \eta_2 = \beta^{\prime} + d \gamma^{\prime}$.   
Then 
$d (\gamma - \gamma^{\prime}) = (\beta^{\prime} - \beta_{Z}) = d\kappa$; 
that is $d (\gamma - \gamma^{\prime} - \kappa) = 0$,
with $\kappa \in \Omega_{\Z}^{n-d-1} (\norm{N \covering})$ integral 
and all forms $\gamma, \gamma^{\prime}, \kappa$ having support in 
$\norm{N \covering} - \norm{N \covering_{X - \abs{Z}}}$. 
But since $H^{n-d-1} (X, X - \abs{Z}) \cong H_{d+1} (\abs{Z}) = 0$, 
we have $\gamma - \gamma^{\prime} - \kappa = d \tau$, for 
$\tau \in \Omega^{n-d-2} (\norm{N \covering})$ also with support 
in $\norm{N \covering} - \norm{N \covering_{X - \abs{Z}}}$. 
Thus we have from \eqref{lambdaequ} 
\begin{equation*}
\begin{aligned}
\gamma^{\prime} + F_1 &= \gamma + F_1 - (\kappa + d \tau)\ , \\
\Lambda_{Z}^{\prime} &= \Lambda_{Z} - (\kappa + d \tau)\ . 
\end{aligned}
\end{equation*}
By Theorem \ref{simpldeligne3}, this shows that the equivalence 
class $[\Lambda_{Z}]$ is well--defined and we get 
\begin{equation}
  \label{lambdaequ3}
[\Lambda_{Z}^{\prime}] = [\Lambda_{Z}] \in H_{\D}^{n-d} (X, \Z)\ . 
\end{equation} 
This proves the theorem for the covering $\covering$ by 
open stars of the vertices of $K$. For an arbitrary good 
covering $\covering^{\prime}$ finer than $\covering$, 
we just use the image of $\Lambda_{Z}$ by the natural map 
of simplicial deRham complexes induced by the map 
$N\covering^{\prime} \to N\covering$. 
It is straightforward to check that this pull--back 
agrees with the construction of $\Lambda_{Z}$ as above 
relative to the refinement $\covering^{\prime}$. 
In particular, the pull--back of the simplicial form $F_1$ 
to $\norm{N\covering^{\prime}}$ coincides with the simplicial 
form $F_1^{\prime}$ defined with respect to $\covering^{\prime}$. 

\noindent
The properties of $\Lambda_{Z}$ stated in $(1)$ to $(4)$ are 
clear from the construction. 
\end{proof}

\begin{remarks}
  \label{lambdaF} 
	
	\paritem{$(1)$} 
The invariance of $[\Lambda_{Z}]$ under refinements applies in 
particular to the covering $\covering^{\prime}$ 
by the open stars of the vertices of a subdivision $K^{\prime}$ of $K$, 
which is a refinement of the covering 
by the open stars of the vertices of $K$. 
	
	\paritem{$(2)$} 
>From the preceding proof, we have 
$F_1 \vert_{\norm{N\covering_{Z}}} = 0$ and therefore 
\begin{equation}
  \label{lambdaequ4}
\Lambda_{Z}\vert_{\norm{N\covering_{Z}}} = 
\gamma \vert_{\norm{N\covering_{Z}}}\ . 
\end{equation}

\end{remarks}

\begin{definition}
  \label{abelgerbetwo}  
  The Deligne cohomology class   
\begin{equation}
   \label{lambdaz} 
   [\Lambda_{Z}] \in H_{\Deligne}^{n-d}(X, \numbersZ)   
\end{equation}
we shall call the \textit{Abel gerbe} associated to the cycle $Z$. 
\end{definition}

\begin{remark}
  \label{submanifold}
In particular, if $M = M^{d} \subset X$ is a closed oriented 
submanifold, we choose a triangulation $K$ of $X$, such that 
$M^d \approx \abs{L}$, where $L$ is a subcomplex of $K$. 
Then there is a canonical simplicial cycle 
$Z_{M} \in C_d (L) \subseteq C_d (K)$, 
such that $\abs{Z_{M}} = \abs{L} \subseteq \abs{K}$ and $Z_{M}$ 
represents the fundamental class $[M] \in H_d (M) \cong \Z$. 
Viewed as a cycle on $X \approx \abs{K}$, $Z_{M} \in C_d (K)$ 
represents the image of $[M]$ under the homomorphism 
$H_d (M) \cong \Z \to H_d (X)$, also denoted by $[M]$; 
that is, we have $[Z_{M}] = [M] \in H_d (X)$. 
Then we put $\Lambda_{M} = \Lambda_{Z_{M}}$, so that the Abel gerbe 
\begin{equation}
   \label{lambdam} 
   [\Lambda_{M}] \in H_{\Deligne}^{n-d}(X, \numbersZ) 
\end{equation}
is well--defined, with $d\Lambda_{M}=\epsilon^*\eta_{M} -\beta_{M}$ 
having the obvious meaning, namely $\eta_{M} = \eta_{Z_{M}}$ and  
$\beta_{M} = \beta_{Z_{M}}$.   

\end{remark}


\section{Linear equivalence of cycles and Abel's Theorem}
\label{secfive}

For $X$ a Riemann surface and $Z = \Divisor$ a divisor as in 
Section \ref{sectwo}, the Abel gerbe is the associated holomorphic 
line bundle $\lineBundle(\Divisor)$ with the complex connection 
given by the holomorphic structure. In this case, by 
Lemma \ref{parallelsection} $Z$ has a meromorphic solution, 
i.e. it is linearly equivalent to zero, 
if and only if $[\Lambda_{Z}] = 0$ in $H^{2}_{\Deligne} (X, \Z)$. 
Motivated by this, we introduce the following definition of 
linear equivalence for cycles.  

\begin{definition}
  \label{LinEq1}
  Two cycles $Z_1, Z_2 \in C_d(K)$ are called \textit{linearly equivalent}
  if 
\begin{equation}
  \label{LinEq2}  
  [\Lambda_{Z_1 - Z_2}] = [\Lambda_{Z_1}]-[\Lambda_{Z_2}] =0 \in 
  H_{\Deligne}^{n-d}(X, \numbersZ)\ .
\end{equation}
\end{definition}

\begin{remark}
  If $[\Lambda_Z] = 0$ then in particular $\eta_Z = 0$ and $[\beta_Z] = 0$
  in $H^{n-d}(X,\numbersZ)$, that is, $Z$ is homologous to zero.
\end{remark}

\begin{theorem}[Abel's Theorem]
  \label{abelone}
  Let $Z=\partial\Gamma$, $\Gamma\in C_{d+1}(K)$. Then $Z$ is linearly
  equivalent to zero, if and only if 
\begin{equation*}   
  \int_\Gamma \theta\in\numbersZ\ , 
\end{equation*}
  for all harmonic $\theta \in {\H}^{d+1}(X, \Z)$ with integral periods.
\end{theorem}
For the proof, we again solve the distributional equation 
\eqref{poisson1} with $\eta_Z = 0$: 
\begin{equation}
  \label{deltagamma}   
  \Delta H_{Z} = -\delta_{Z} = - d \delta_{\Gamma}\ , 
\end{equation}
where $\delta_{\Gamma} (\psi) = \int_{\Gamma} \psi, 
\psi \in \Omega^{d + 1} (X)$. 
Hence for $F_{\partial\Gamma} = F_{Z} = \ast~d \ast H_{Z}$ as before, 
we get 
\begin{equation*}
\Delta H_{Z} =  d F_{\partial\Gamma} = - d \delta_{\Gamma}\ , 
\end{equation*}
and
\begin{equation*}
d ( F_{\partial\Gamma} + \delta_{\Gamma} ) = 0\ , 
\end{equation*}
so that by de Rham--Hodge theory for currents~\cite{DRh} 
\begin{equation}
  \label{alphagamma}   
  F_{\partial\Gamma} + \delta_{\Gamma} = \alpha_{\Gamma} + d T\ , 
\end{equation}
for a \textit{harmonic} form $\alpha_{\Gamma} \in \mathcal{H}^{n - d -1} (X)$ 
and a $(n-d-2)$--current $T\in \Omega^{d + 2}(X)^{\prime}$. 
Note that $\alpha_{\Gamma}$ is smooth by elliptic regularity. 

We shall first prove the following theorem. 
\begin{theorem}
\label{abeltwo}
  For $Z = \partial\Gamma$, the simplicial gerbe $\Lambda_{Z}$ and the harmonic 
  form $\alpha_{\Gamma}$ have the following properties: 
  
  \paritem{$(1)$} 
  As simplicial forms, we have 
\begin{equation}
  \label{lambdaalpha}   
  \Lambda_{Z} = \Lambda_{\partial \Gamma} \equiv 
  \epsilon^{\ast} \alpha_{\Gamma} \mod \left( \Omega^{n-d-1}_{\Z}(\norm{N \covering}) + 
  d \Omega^{n-d-2} (\norm{N \covering}) \right)\ ; 
\end{equation}  
  that is, the simplicial form $\Lambda_{\partial\Gamma}$ 
  is given by the global harmonic form $\alpha_{\Gamma}$. 
  
  \paritem{$(2)$} 
  There exists an integral form $\kappa \in \Omega^{n-d-1}_{\Z}(\norm{N \covering})$ 
  with support in a regular neighborhood $V_{\Gamma}$ of $\abs{\Gamma}$, 
  such that for all harmonic $(d+1)$--forms with integral periods 
  $\theta \in \H^{d+1} (X, \Z)$, we have   
\begin{equation}
  \label{lambdakappa}
  \int_{[X]} (\Lambda_{\partial\Gamma} + \kappa) \wedge \epsilon^{\ast} \theta \equiv 
  \int_{\Gamma} \theta \mod \Z \ . 
\end{equation} 
  
  \paritem{$(3)$} 
  If $Z = \partial\Gamma = \partial\Gamma^{\prime}$, then 
  $\zeta = \alpha_{\Gamma^{\prime}} - \alpha_{\Gamma} \in 
  \H^{n-d-1} (X,\Z)$; that is, $\zeta$ is a harmonic form 
  with integral periods. 
  Hence, $[\alpha_{\Gamma}]$ is well-defined in the Picard 
  torus $\Pic^{n-d-1}(X) = \H^{n-d-1} (X) / \H^{n-d-1} (X,\Z)$ 
  in \eqref{picardone}.
  
\end{theorem}

\begin{proof}
First notice that since $F_{Z} = \ast~d \ast H_{Z}$ and $\theta$ is harmonic, 
we get from \eqref{alphagamma}
\begin{equation}
  \label{alphatheta}   
  \int_{X} \alpha_{\Gamma} \wedge \theta = 
  \innerp{ F_{Z} + \delta_{\Gamma}, \theta } = 
  \innerp{  \delta_{\Gamma} , \theta } = \int_{\Gamma} \theta\ . 
\end{equation}
This shows that $(1)$ and $(2)$ are equivalent. 

For the proof of $(2)$, we let $V_{\Gamma} = 
\bigcup \{ U_i \in \covering_{\Gamma}\}$, where  
$\covering_{\Gamma}$ is the set of open sets $U_i \in \covering$ 
intersecting $\Gamma$, so that $V_{\Gamma}$ is a regular neighborhood
of $\abs{\Gamma}$. Since formula \eqref{lambdakappa} is additive in $\Gamma$, 
we can without loss of generality assume that $\Gamma$ consists of a 
single simplex and that $V_{\Gamma}$ is contractible. Therefore we can
assume that $\theta \vert_{V_{\Gamma}} = d \nu$ for some 
$\nu \in \Omega^d (V_{\Gamma})$. 
>From the formulas for integration of simplicial forms 
(cf. Dupont--Kamber~\cite{DK2} and Dupont--Ljungmann~\cite{DL} ), 
together with the construction of $\Lambda_{Z}$ in the proof of 
Theorem \ref{abelgerbeone}, we now get
\begin{equation}
\label{alphatheta2}
  \begin{aligned}
  \innerp{  F_{Z} + \delta_{\Gamma}, \theta } &=  
  \int_{[X - V_{\Gamma}]} F_{Z} ~\wedge \theta + 
  \innerp{  F_{Z}\vert_{\overline{V}_{\Gamma}} + \delta_{\Gamma}, d \nu }     \\
&=  \int_{[X - V_{\Gamma}]} \Lambda_{Z} \wedge \veps^{\ast} \theta - 
  \int_{[\partial\overline{V}_{\Gamma}]} \Lambda_{Z} \wedge \veps^{\ast} \nu   \\
&= \int_{[X - V_{\Gamma}]} \Lambda_{Z} \wedge \veps^{\ast} \theta + 
   \int_{[\overline{V}_{\Gamma}]} d \kappa \wedge \veps^{\ast} \nu + 
   \int_{[\overline{V}_{\Gamma}]} \Lambda_{Z} \wedge \veps^{\ast} \theta  \\
&= \int_{[X]} \Lambda_{Z} \wedge \veps^{\ast} \theta + 
   \int_{[\overline{V}_{\Gamma}]} \kappa \wedge \veps^{\ast} \theta + 
   \int_{[\partial\overline{V}_{\Gamma}]} \kappa \wedge \veps^{\ast} \nu \\
&= \int_{[X]} (\Lambda_{Z} + \kappa) \wedge \veps^{\ast} \theta\ . 
  \end{aligned}
\end{equation}
Here we have used that, since $Z = \partial\Gamma \sim 0$, we have $\eta_{Z} = 0$ 
and $\beta_{Z} = - d \kappa$ and hence $d \Lambda_{Z} = d \kappa$ for some 
integral simplicial form $\kappa$ with support in $V_{\Gamma}$. 
We used also the simplicial Stokes' theorem~\cite{DL} to see that 
$\int_{[\overline{V}_{\Gamma}]} d (\kappa \wedge \veps^{\ast} \nu) 
=  \int_{[\partial\overline{V}_{\Gamma}]} \kappa \wedge \veps^{\ast}
\nu = 0$, since $\kappa$ vanishes on $\partial\overline{V}_{\Gamma}$. 
Equations \eqref{alphatheta} and \eqref{alphatheta2} now prove $(2)$. 

For the proof of $(3)$, let $Z = \partial\Gamma = \partial\Gamma^{\prime}$. 
Then 
$\partial (\Gamma^{\prime} - \Gamma) = 0$ and $Z^{\prime} = 
\Gamma^{\prime} - \Gamma$ is an integral $(d+1)$--cycle. 
Equation \eqref{alphagamma} implies
$
\delta_{Z'} = \delta_{\Gamma'} - \delta_{\Gamma} = 
(\alpha_{\Gamma'} - \alpha_{\Gamma}) + d( T'-T).  
$
So $\zeta = \alpha_{\Gamma^{\prime}} - \alpha_{\Gamma}$ satifies 
$\delta_{Z'} = \zeta + d(T'-T)$. 
Since $Z'$ is an integral cycle, $\zeta$ must be an integral harmonic 
form $\zeta \in \H^{n-d-1} (X,\Z)$. 
\end{proof}

\medbreak
Abel's Theorem \ref{abelone} is now a consequence of the following 
Corollary to Theorem \ref{abeltwo}. 

\begin{corollary}
\label{abelthree}
  For $Z = \partial \Gamma$ as above, the following statements are 
  equivalent:
  
  \paritem{$(1)$} 
  $[\Lambda_{Z}] = 0$ in $H_{\Deligne}^{n-d}(X, \numbersZ)$; 
  
  \paritem{$(2)$} 
  For all harmonic $(d+1)$--forms $\theta$ with integral periods; 
  that is, $\theta \in \mathcal{H}^{d+1} (X, \numbersZ)$, we have   
\begin{equation}
  \label{gammatheta} 
  \int_{\Gamma} \theta \in \numbersZ \ . 
\end{equation}

  \paritem{$(3)$} 
  There exists $\Gamma_{0}$ with $\partial \Gamma_{0} = Z$, such that 
\begin{equation*}
  F_{\partial\Gamma_{0}} + \delta_{\Gamma_{0}} = d T_{0}\ , 
\end{equation*}
  where $F_Z = F_{\partial \Gamma_{0}}$ is given as before. 
  By \eqref{alphagamma}, we have $\alpha_{\Gamma_{0}} = 0$.
\end{corollary}

\begin{proof}
By Theorem \ref{abeltwo} $[\Lambda_{Z}]$ is represented in 
$H^{n-d-1} (X, \numbersR)$ by the harmonic form $\alpha_{\Gamma}$. 
Hence (1) and (2) are equivalent to $\alpha_{\Gamma} \equiv 0 
\mod H^{n-d-1} (X, \numbersZ)$. 
>From Theorem \eqref{abeltwo} (3), we know that $[\alpha_{\Gamma}] = 0 \in 
\H^{n-d-1} (X) / \H^{n-d-1} (X,\Z)$.  
By changing $\Gamma$ by a cycle, we 
can make $\alpha_{\Gamma} = 0$. This proves that (3) is equivalent 
to (1) and (2). 
\end{proof}

\begin{remark}
\label{meromorphic}
Notice that $F = F_{\partial\Gamma} \vert_{X - \abs{Z}}$ is now 
harmonic in the stronger sense that $dF = 0$ and $d^{\ast} F = 0$
by \eqref{poisson3}. Thus $F_{\partial\Gamma}$ is analogous to a 
meromorphic solution in the classical Abel Theorem. 

\end{remark}


\section{Moduli spaces}
\label{secsix}
In this section, we need to enlarge the chain complex 
$C_{\ast} (K)$ relative to a smooth triangulation of $X$, 
that is $X \approx \abs{K}$, 
which was introduced at the beginning of Section \ref{secfour}. 
Therefore we look at the limit complex 
\begin{equation}
  \label{limittriang}
  C_{\ast} (X) = \lim_{\underset{K} \lra} ~C_{\ast} (K) 
  \subset \mathcal{S}_{\ast} (X)\ , 
\end{equation}
taking into account the inclusions of chain complexes 
$C_{\ast} (K) \subseteq C_{\ast} (K^{\prime})$ where $K^{\prime}$ 
corresponds to a subtriangulation of the triangulation coming from $K$. 
Obviously, we can view $C_{\ast} (X)$ as a subcomplex of the singular 
complex $\mathcal{S}_{\ast} (X)$ of $X$. Then 
$C_{\ast} (K) \subseteq C_{\ast} (K^{\prime}) \subset \mathcal{S}_{\ast} (X)$
induce isomorphisms in homology, so that we have canonical isomorphisms 
\begin{equation*}
H_{\ast} (C_{\ast}(X)) \cong \lim_{\underset{K} \lra} ~
H_{\ast} (K) \cong H_{\ast} (X)\ . 
\end {equation*}
The construction of the Abel gerbe $Z \to [\Lambda_{Z}]$ passes 
to the limit \eqref{limittriang} and defines a homomorphism 
$Z_d (X) \to H_{\D}^{n-d} (X, \Z)$. 
This follows from the proof of Theorem \ref{abelgerbeone} and the 
fact that the covering $\covering^{\prime}$ given by the open stars 
of the vertices of a subdivision $K^{\prime}$ of $K$ is a refinement 
of the covering $\covering$ given by the open stars of the vertices of $K$. 

\medbreak
As the construction of the Abel gerbe in Section \ref{secfour} 
involves deRham--Hodge theory on the compact oriented Riemannian 
manifold $X$, we need now to better understand the terms in diagram 
\eqref{DeligneDiagram} for the Deligne cohomology in view 
of the deRham--Hodge decomposition of forms on $X$: 
\begin{equation}
   \label{DRH} 
\Omega^{\ell} (X) \cong \mathcal{H}^{\ell} (X) \oplus d \Omega^{\ell -1}(X) 
\oplus d^{\ast} \Omega^{\ell +1}(X)\ . 
\end{equation}
We recall that $\mathcal{H}^{\ell} (X, \mathbb{Z} ) \subset \mathcal{H}^{\ell} (X)$ 
denotes the harmonic forms, respectively the integral lattice of harmonic forms. 
Further, the sum decompositions in \eqref{DRH} and the following formulas 
are orthogonal.
Thus the deRham--Hodge decomposition \eqref{DRH} implies that
\begin{equation}
   \label{DRH1} 
H_\D^{\ell+1} (X) \cong \mathcal{H}^{\ell} (X) \oplus d^{\ast} \Omega^{\ell +1}(X) \quad , 
\quad \Omega^\ell (X) / \Omega_{\rm cl}^{\ell} (X) \cong d^{\ast} \Omega^{\ell +1}(X) 
\end{equation}
and also 
\begin{equation}
   \label{DRH2} 
\Omega_{\rm cl}^{\ell} (X ,\mathbb{Z}) \cong 
\H^{\ell} (X, \Z) \oplus d d^{\ast} \Omega^{\ell}(X)\ . 
\end{equation}
This implies
\begin{equation}
   \label{DRH3} 
j_{\ast} H^{\ell} (X,\Z) \cong \Omega_{\rm cl}^{\ell} (X ,\Z) / 
d d^{\ast} \Omega^{\ell}(X) \cong \H^{\ell} (X, \Z)\ , 
\end{equation}
as well as 
\begin{equation}
   \label{DRH4} 
\Omega^{\ell} (X)  / \Omega_{\rm cl}^{\ell} (X ,\mathbb{Z}) \cong 
\H^{\ell} (X) / \H^{\ell} (X, \Z) \oplus d^{\ast} \Omega^{\ell +1}(X)\ . 
\end{equation}

\begin{remarks}
  \label{delignerem2} 

This has the following consequences for the diagram \eqref{DeligneDiagram}: 

	
	\paritem{$(1)$} 
	By \eqref{DRH1}, the right arrow in the second exact row is of the form
\begin{equation}
  \label{projection}
	H_\D^{\ell+1} (X) \cong \mathcal{H}^{\ell} (X) \oplus d^{\ast} \Omega^{\ell +1}(X) 
	\rightarrow \Omega^\ell (X) / \Omega_{\rm cl}^{\ell} (X) \cong d^{\ast} \Omega^{\ell+1}(X)\ , 
\end{equation} 
and is given by orthogonal projection to the second summand. Here the infinite 
dimensional part $d^{\ast} \Omega^{\ell +1}(X)$ consists of topologically 
trivial gerbes of the form $\omega_{0} = d^{\ast} \alpha$ whose curvature 
$d \omega_{0} = d d^{\ast} \alpha$ uniquely determines 
$\omega_{0} = d^{\ast} \alpha$. 
	
	\paritem{$(2)$} 
	Using \eqref{DRH1}, \eqref{DRH3}, the kernel of $\iota_{\ast}$ are the 
	harmonic forms $\mathcal{H}^{\ell} (X, \Z)$ with integral periods. 
	Thus the image of $\iota_{\ast}$ contains the torus 
	\begin{equation}
    \label{torus}
    \vcenter{
  \xymatrix{
  H^\ell (X,\R) / j_{\ast} H^\ell (X,\Z) \ar[d]^{\subseteq}_{\overline{\rho}_{\ast}} 
  \ar[r]^-{\cong}  & \mathcal{H}^{\ell} (X) / \mathcal{H}^{\ell} (X, \Z) 
  \ar[d]^{\subseteq}_{\overline{\iota}_{\ast}}  \\
  H^\ell (X, \R / \Z)  \ar[r]^{\subseteq}  & H_\D^{\ell+1}(X,\Z) 
  }
}
  \end{equation}
	of topologically trivial \textit{flat} $\ell$--gerbes. 
  In our motivating situation in Section \ref{sectwo}, where 
  $\ell=1$, this torus corresponds to the Picard variety of 
  topologically trivial holomorphic line bundles. 
  We will refer to it as the \textit{Picard torus} and write 
  \begin{equation}
    \label{picardone}
  \Pic^{\ell} (X) = \mathcal{H}^{\ell} (X) / \mathcal{H}^{\ell} (X, \Z)\ . 
  \end{equation}
  
  Note that from \eqref{DeligneDiagram} and Remark \ref{delignerem} $(3)$, 
	the Picard torus in \eqref{torus} differs from the moduli space 
	$H^\ell (X, \R / \Z)$ of flat $\ell$--gerbes by the torsion subgroup 
	of $H^{\ell +1} (X, \Z)$. 
	This is encoded in diagram \eqref{DeligneDiagram} by the left exact 
	column; that is, the Bockstein exact sequence. In fact, the torus on 
	the left side of \eqref{torus} is exactly the kernel of the Bockstein
	boundary map $\beta_{\ast}$ and the image of $\beta_{\ast}$ is the 
	torsion subgroup of $H^{\ell +1} (X, \Z)$. 
 	
  \paritem{$(3)$} 
	It follows from \eqref{DeligneES2} and \eqref{DRH4} that the Deligne 
	cohomology is given by an exact sequence 
	(i.e. the middle exact column in \eqref{DeligneDiagram})
\begin{equation}
  \label{DeligneES3}
    0 \ra \Pic^{\ell} (X) \oplus d^{\ast} \Omega^{\ell +1}(X)  
    \xra[\subseteq]{\overline{\iota}_{\ast}} 
    H_\D^{\ell+1}(X,\Z) \xra{c} H^{\ell +1}(X,\Z) \ra 0\ , 
\end{equation}
  since $\overline{\iota}_{\ast}$ is injective on $d^{\ast} \Omega^{\ell +1}(X)$ by exactness 
  of the third column of \eqref{DeligneDiagram}. 

  \paritem{$(4)$} ~Harmonic Deligne cohomology: 
  If we pull back the exact sequence \eqref{DeligneES} along the inclusion 
  ${\H}^{\ell +1} (X, \Z) \subset \Omega_{\rm cl}^{\ell+1}(X,\Z)$, 
  we obtain the \textit{harmonic} Deligne cohomology ${\H}_{\D}^{\ell+1}(X, \Z)$ 
  of $\ell$--gerbes with \textit{harmonic curvature}: 
\begin{equation}
   \label{redDeligne}
   \vcenter{
   \xymatrix{
 0 \ar[r] & H^\ell(X,\R / \Z)  \ar[r] \ar[d]^{\id} & {\H}_{\D}^{\ell+1}(X, \Z) 
 \ar[r]^-{d_*} \ar[d]^{\subset} & {\H}^{\ell +1} (X, \Z) \ar[r] 
 \ar[d]^{\subset} &0   \\
0 \ar[r] & H^\ell(X,\R / \Z) \ar[r]  & H_{\D}^{\ell+1}(X,\Z) 
\ar[r]^-{d_*} & \Omega_{\rm cl}^{\ell+1}(X,\Z) \ar[r] & 0  
  }
}
\end{equation}
Then the exact sequence \eqref{DeligneES3} becomes 
\begin{equation}
  \label{DeligneES4}
    0 \ra \Pic^{\ell} (X) \xra[\subseteq]{\overline{\iota}_{\ast}} 
    \H_\D^{\ell+1}(X,\Z) \xra{c} H^{\ell +1}(X,\Z) \ra 0\ , 
\end{equation} 
For these reasons, we call these gerbes \textit{harmonic} gerbes 
and ${\H}_{\D}^{\ell+1}(X, \Z)$ the \textit{harmonic} Deligne 
cohomology. 


\end{remarks}

\subsection{The Picard torus and the Picard map}
\label{subsecsixone}

  From the construction of the Abel gerbe in Theorem \ref{abelgerbeone} 
  and the definition of linear equivalence in Definition \ref{LinEq1},
  we have an injection of abelian groups
  \begin{equation}
  \label{lambdabar}
  \overline{\Lambda} \colon \M_d (X) \colon = Z_d (X)/\{\text{lin.\ equiv.}\} 
  \subseteq \H_{\D}^{n-d}(X,\numbersZ)\ , 
  \end{equation}
  where the inclusion $\overline{\Lambda}$ is induced by 
  $Z \mapsto [\Lambda_{Z}]$.   

  For the boundaries $B_d (X) \subset Z_d (X)$, 
  we have the following inclusion from Theorem \ref{abeltwo} (1), (3):
  \begin{equation}
    \label{alphabar}
  \overline{\alpha} \colon \M_d^{\circ} (X) \colon = B_d(X)/\{\text{lin.\ equiv.}\} 
  \subseteq \Pic^{n-d-1} (X)\ , 
  \end{equation}
  where the inclusion is given by $Z = \partial \Gamma \mapsto 
  \overline{\alpha}_{\Gamma} = [\alpha_{\Gamma}]$. 
  Thus $\M_{d} (X)$, respectively $\M_{d}^{\circ} (X)$, is the 
  \textit{moduli space of Abel gerbes}, respectively the 
  \textit{moduli space of topologically trivial Abel gerbes}. 
  From Theorem \ref{abeltwo} (1) and \eqref{torus} we have 
  the following Cartesian diagram; that is, a pull--back diagram: 
  \begin{equation}
  \label{cartesian}
  \vcenter{
    \xymatrix{
      \M_d (X) \ar@{>->}[r]^-{\overline{\Lambda}}_-{\subseteq}  & 
      \H_{\Deligne}^{n-d}(X,\numbersZ) \\ 
      \M_{d}^{\circ} (X) \ar@{>->}[r]^-{\overline{\alpha}}_-{\subseteq}  
      \ar[u]^{\subseteq} & \Pic^{n-d-1} (X) 
      \ar[u]^{\subseteq}_{\overline{\iota}_{\ast}}  
    }
  } 
  \end{equation} 

Recall that by construction, the image of $\M_{d} (X)$ in 
  $H_{\Deligne}^{n-d}(X,\numbersZ)$ is contained in the group 
  of gerbes whose curvature is harmonic with integral periods; 
  that is, in the harmonic Deligne cohomology 
  $\H_{\Deligne}^{n-d}(X,\numbersZ)$ (cf. \eqref{redDeligne}). 
  In contrast, $\M_{d}^{\circ} (X)$ is exactly the part of 
  $\M_{d} (X)$ which maps into the Picard torus \eqref{torus}, 
  \eqref{picardone}, namely $\Pic^{n-d-1} (X)$; 
  that is, it consists of flat, topologically trivial gerbes. 
  We call $\overline{\alpha}$ the \textit{Picard} map. 
    
  From \eqref{DeligneES4} and the fact that the characteristic 
  class of the Abel gerbe $[\Lambda_Z]$ is the Poincar\'e dual 
  $[\beta_Z]$ of $[Z]$, it follows that we have canonical isomorphisms
  \begin{equation}
    \label{Pdual}
  \M_d (X) / \M_{d}^{\circ} (X) \cong H_{d} (X, \numbersZ) 
  \xra[\cong]{PD} H^{n-d}(X,\numbersZ)\ , 
  \end{equation}
  the second being Poincar\'e duality, induced by the characteristic class. 
  We will see in Proposition \ref{deftwo} that $\M_d^{\circ} (X)$ is 
  connected in $\Pic^{n-d-1}(X)$. So if $\M_d^{\circ} (X) \cong \Pic^{n-d-1}(X)$, 
  then $\M_d (X) \cong \H_{\D}^{n-d} (X, \Z)$, the harmonic Deligne cohomology 
  of classes with harmonic curvature. 

Thus, we need to understand the image of $\mathcal{M}_{d}^{\circ}$ 
in the Picard torus $\Pic^{n-d-1} (X)$ of topologically trivial flat gerbes. 

\begin{remark}
  \label{torsion} ~Torsion classes (cf. Remark \ref{delignerem2} (2)): 
Suppose that the Abel gerbe $\Lambda_{Z}$ is flat; 
that is $\eta_{Z} = 0$, so that 
$[\Lambda_{Z}] \in H^{n-d-1} (X, \numbersRmodZ)$. By diagram 
\eqref{DeligneDiagram}, the characteristic class 
$[\beta_{Z}] \in H^{n-d} (X, \Z)$ is given by 
$[\beta_{Z}] = \beta_{\ast} [\Lambda_{Z}]$, where $\beta_{\ast}$ 
is the Bockstein homomorphism. Thus $\beta_{Z}$ is a torsion class, 
say $m \cdot [\beta_{Z}] = 0$ for some $m \in \N^{+}$. 
By Poincar\'e duality, we have also $m \cdot Z = \partial\Gamma$ and
so $m \cdot Z$ determines an element in $\M_{d}^{\circ} (X)$. 
Finally the Bockstein formula implies 
$m \cdot \beta_{\ast} [\Lambda_{Z}] = 0$; 
that is, $m \cdot [\Lambda_{Z}] = [\alpha_{\Gamma}]$ 
takes value in the Picard torus $\Pic^{n-d-1} (X)$. 

\end{remark}


\subsection{The Jacobi torus and the Abel--Jacobi map}
\label{subsecsixtwo}

First, we observe that $\alpha \mapsto \int_{X} \alpha \wedge$ induces by 
Poincar\'e duality a canonical isomorphism 
$\varphi \colon \H^{n-d-1}(X) \cong \H^{d+1}(X)^{\ast}$. It further induces 
an isomorphism of abelian tori of (real) dimension $\dim H^{d+1} (X. \R)$~:
\begin{equation}
  \label{abeljacobi1}  
  \overline{\varphi} \colon 
  \mathcal{H}^{n-d-1} (X) / \mathcal{H}^{n-d-1} (X, \Z) \cong  
  \Hom (\mathcal{H}^{d+1} (X, \Z), \numbersRmodZ)\ .   
\end{equation}
This is valid for $d=0,\dots,n-1$. 
The torus on the right hand side of \eqref{abeljacobi1} 
corresponds classically to the Jacobi variety of a Riemann surface, 
where $n=2,~d=0$. 
We shall call it the 
\textit{Jacobi torus} and denote it by $\J^{d+1}(X)$. 
We now recall formula \eqref{alphatheta}; that is, 
\begin{equation*}
  \int_{X} \alpha_{\Gamma} \wedge \theta = \int_{\Gamma} \theta, ~
  \theta \in \H^{d+1}(X)\ . 
\end{equation*}
Combining \eqref{alphatheta} with \eqref{abeljacobi1}, we obtain a 
commutative diagram 
\begin{equation}
  \label{abeljacobi2}   
  \vcenter{
    \xymatrix{ 
      &  \Pic^{n-d-1}(X) \ar[dd]^{\overline{\varphi}}_{\cong}   \\
      \M_d^{\circ} (X) \ar@{>->}[ur]^{\overline{\alpha}}_{\subseteq} 
      \ar@{>->}[dr]^{\overline{J}}_{\subseteq}  \\
      &  \J^{d+1}(X) 
    }
  }
\end{equation} 
where $\overline{J}$ is induced by the functional 
\begin{equation}
  \label{abeljacobi}
Z = \partial\Gamma \mapsto 
J_{\partial\Gamma} (\theta) = \int_{\Gamma} \theta, 
~\theta \in \mathcal{H}^{d+1} (X, \Z)\ . 
\end{equation}
Note that $\overline{J}$ is well-defined and injective by 
Abel's Theorem \ref{abelone}, plus the fact that 
$J_{\partial\Gamma}(\theta)$ has integral values if $\Gamma$ is a 
cycle; that is $\partial\Gamma = 0$. 
We shall call $\overline{J}$ the \textit{Abel--Jacobi map}. 
Thus we may just as well use the map $\overline{J}$ to investigate 
the image of $\M_{d}^{\circ} (X)$. 
The Abel--Jacobi map $\overline{J}$ is given in terms of period integrals and 
therefore is more explicit than the Picard map $\overline{\alpha}$,  
which is determined by the solution of a Laplace--Poisson equation. 
Therefore it is in general easier to deal with and more effective 
in explicit calculations, as we shall see. 

\begin{remarks}
  \label{intermediate} ~Intermediate Jacobians: 

\paritem{$(1)$} 
For a K\"ahler manifold, our definition of the Jacobians agrees 
with the tori underlying the complex \textit{intermediate} Jacobians in odd 
degrees, which are related to \textit{holomorphic} Deligne cohomology 
(cf. Griffiths--Harris \cite{GH}, Ch. 2.6, Dupont--Hain-Zucker \cite{DHZ} 
and also Clemens~\cite{Cle}, Harvey--Lawson~\cite{HL3} and the references 
given there). 
For divisors on algebraic manifolds of complex dimension greater than one, 
it is not clear how our version of Abel's theorem is related to 
the version by Griffiths \cite{Gr}. 

\paritem{$(2)$} 
We also remark that for $\dim X = n = 4k+2, ~d = 2k, ~k \geq 0$,
the Picard and the Jacobi tori in degree $n-d-1 = d+1 = 2k+1$ 
carry a canonical \textit{complex} and \textit{symplectic} structure, 
compatible with the isomorphism 
\begin{equation*}   
\Pic^{2k+1}(X) \xra[\cong]{\overline{\varphi}} \J^{2k+1}(X)\ . 
\end{equation*} 
The former is induced by the Hodge $\ast$--operator on $\H^{2k+1}(X)$ 
and the latter is defined by the pairing 
$\innerp{ \alpha, \beta } = \int_X \alpha \wedge \beta$ 
on $\H^{2k+1}(X)$. 
\end{remarks}


\subsection{Deformations}
\label{subsecsixthree}
We now consider `deformations' of Abel gerbes as follows: 

\begin{definition}
  \label{defone} 

\noindent
\paritem{$(1)$}
  A \textit{regular $(d+1)$--simplex} in $X$ is a smooth embedding 
  $\Gamma \colon \simplex^{d+1} \to X$ of the standard simplex 
  $\simplex^{d+1}$ 
  (or rather an open neighborhood in the hyperplane 
  $\sum_{i=0}^{d+1} t_i = 1$). 
  Note that any simplex in a triangulation $K$ of $X$ can be 
  parametrized as a regular simplex.

\paritem{$(2)$}
  A \textit{deformation} of a cycle $Z = \partial\Gamma \in B_d (X)$ 
  for a triangulation $K$ of $X$ is a family of cycles 
  $Z_{r} = \partial\Gamma_{r} \in B_d (X),~r \in [0,1]$, 
  for some subdivisions $K_{r}$ of $K$, 
  such that $\Gamma_1 = \Gamma$, ~
  each simplex of $\Gamma_r$ is regular and 
  $r \mapsto \overline{J} (Z_{r}) = J_{\partial\Gamma_{r}} = 
  \int_{\Gamma_{r}} (\cdot) \in \J^{d+1} (X)$ is a smooth curve. 
      
\end{definition}
The following deformation techniques are used repeatedly in what 
follows and we state them in a separate Lemma.  

\begin{lemma}
  \label{variation}
Let $Z = \partial \Gamma$, for $\Gamma$ any $(d+1)$--chain 
in $C_{d+1} (K)$ consisting of regular simplices, 
where $K$ is an arbitrary triangulation of $X$. Then there is a 
deformation $Z_r = \partial\Gamma_r, ~r \in [0,1]$ of $Z$ satisfying: 

  \paritem{$(1)$} 
  $\Gamma_r \in C_{d+1} (K_r), ~\abs{K_r}$ a subdivision of $\abs{K}$.
  
  \paritem{$(2)$} 
  For $\alpha_r = \alpha_{\Gamma_r}$, 
  the map $]0,1] \ra \Pic^{n-d-1} (X)$ given by 
  $r \mapsto [\alpha_r]$ is smooth.
  
  \paritem{$(3)$} 
  For $J_r = J_{\partial\Gamma_r}$,  
  the map $]0,1] \ra \J^{d+1} (X)$ given by 
  $r \mapsto \overline{J}_r$ is smooth.
    
  \paritem{$(4)$} 
  $[\alpha_r] \ra 0$ in $\Pic^{n-d-1} (X)$ for $r \downarrow 0$. 
  
  \paritem{$(5)$} 
  $J_r (\theta) =  \int_{\Gamma_r} \theta \ra 0, ~r \downarrow 0$, 
  for all $\theta \in \H^{d+1} (X)$. 

\end{lemma}

\begin{proof} 
 It is clearly enough to take $\Gamma$ to be a regular 
 $(d+1)$--simplex $\Gamma \colon \simplex^{d+1} \to X$ of $K$. 
 Then we simply define 
 $\Gamma_{r} = \Gamma \circ \phi_{r}$, 
with $\phi_r (t_0, \ldots, t_{d-1}, t_d, t_{d+1}) = 
(t_0, \ldots, t_{d-1}, t_d+(1-r) ~t_{d+1}, r ~t_{d+1}), 
 ~t \in \simplex^{d+1}$. 
   Then $\Gamma_r, ~r \in [0,1]$ clearly satisfies $(1)$ and $(2)$. 
 $(2)$ and $(3)$ are equivalent by formula \eqref{alphatheta}. 
 Furthermore by Theorem \ref{abeltwo} and formula \eqref{alphatheta}, 
 conditions $(4)$ and (5) are equivalent and are fulfilled, since 
 $\int_{\Gamma_r} \theta \ra 0, r \downarrow 0$, 
 for all $\theta \in \H^{d+1} (X)$ by construction of $\Gamma_r$. 
\end{proof}


\begin{proposition}
  \label{deftwo}

\paritem{$(1)$} 
  $\M_d^{\circ} (X)$ is connected in the Picard torus, respectively the 
  Jacobi torus.
  
\paritem{$(2)$} 
  The closure $\overline{\M_d^{\circ}} (X)$ in the induced topology is a subtorus 
  of $\Pic^{n-d-1} (X)$. 
  
\paritem{$(3)$} 
For $\overline{\alpha}$, respectively $\overline{J}$ 
to be surjective, it is necessary and sufficient that their image 
contain an open neighborhood of the origin (or an open neighborhood 
of any point in their image). 
  
\end{proposition}

\begin{proof}
 Again let $\partial \Gamma$ for $\Gamma$ any $(d+1)$--chain in $C_{d+1} (K)$, 
 where $K$ is an arbitrary triangulation of $X$. 
 To prove (1), we again take $\Gamma$ to be a regular $(d+1)$--simplex
 and we define as before
 $\Gamma_r (t) = \Gamma (t_0, \dots, t_{d-1}, t_d+(1-r) t_{d+1}, r ~t_{d+1}), 
 ~t \in \simplex^{d+1}$. 
 Then $(1)$ follows from Lemma \ref{variation} and $(2)$ clearly follows from $(1)$. 
 To prove $(3)$, we have only to observe that $\overline{\alpha}$, respectively 
 $\overline{J}$ are homomorphisms of abelian groups. 
 The statement follows from the fact that any open neighborhood of the
 origin in either torus generates the entire torus. 
 Observe that, except for $(1)$, the above deformations can take place
 in the interior of the fundamental domain of ${\H}^{n-d-1} (X)$ 
 relative to the integral lattice ${\H}^{n-d-1} (X, \Z)$. 
\end{proof}


\subsection{The moduli theorem}
\label{subsecsixfour}
In this Section we determine the moduli space of Abel gerbes 
by establishing an inversion theorem for the Abel--Jacobi map. 
Before stating and proving the main Theorem \ref{modulitheorem}, 
we will illustrate the technique involved in some important examples.

\begin{example}
  \label{example1}
  ~The case $n \geq 2, ~d=n-1 \colon$ 

\noindent
  The Jacobi map $\overline{J} \colon \M_{n-1}^{\circ} (X) \to \J^{n} (X)$ 
  is an isomorphism. Therefore, so is the Picard map 
  $\overline{\alpha} \colon \M_{n-1}^{\circ} (X) \to \Pic^{0} (X)$. 
  This is the easiest case, since the Picard and Jacobi tori 
are now in degree $0$, respectively $n$. Thus 
\begin{equation*}
\Pic^{0}(X) \xra[\cong]{\overline{\varphi}} 
\J^{n}(X) \cong \Hom (\mathcal{H}^{n} (X, \Z), \numbersRmodZ) \cong \numbersRmodZ \ ,
\end{equation*}
with the generator of the integral lattice given by $\theta_{0} = \Vol$, 
assuming that the volume is normalized. Taking a deformation $\Gamma_r, 
~r \in [0,1]$ of a regular $n$--simplex as in Lemma \ref{variation}, we get 
$J_{\Gamma_r} (\theta_{0}) =  \int_{\Gamma_r} \Vol > 0$, respectively 
$J_{\Gamma_r} (\theta_{0}) =  \int_{\Gamma_r} \Vol < 0$, if the orientation 
of the regular simplex $\simplex^{n}$ is reversed. Further, we have 
$J_{\Gamma_r} (\theta_{0}) =  \int_{\Gamma_r} \Vol \ra 0, ~r \downarrow 0$. 
Thus the image contains an interval around the origin and so the Jacobi 
map must be an isomorphism $\M_{n-1}^{\circ} (X) \cong \J^{n} (X)$.  

In this case, the Deligne cohomology $H_{\D}^1 (X, \Z)$ 
consists of $0$--gerbes, which are given by 
$f^{0} \in \setC^{0} (\covering, \underline{\R})$, 
such that $\check\delta f^{0} \equiv 0 \mod \Z$, so that 
$f^{0}$ defines a global smooth function 
$\theta \colon X \to \numbersRmodZ \cong U(1)$, 
modulo global functions $\Omega^{0} (X)$. 
Since $\check\delta d f^{0} = d \check\delta f^{0} = 0$, 
the curvature $F_f$ is a closed $1$--form with integral periods, 
determined by $\epsilon^{\ast}( F_f ) = d f^{0}$. 
The characteristic class of $[\theta, f]$ is given by 
$c [\theta, f] = [\check\delta f^{0}]  \in H^1(X, \Z)$; 
that is, the obstruction to lift $\theta$ to a global function 
$f \in \Omega^{0} (X)$. 

\end{example}

\begin{example}
  \label{example2}
  ~The case $n > 2, ~d = n - 2 \colon$ 

\noindent
Here we look at Abel $1$--gerbes associated to submanifolds 
$M^{n-2} \subset X^n$ of codimension $2$ or more generally to 
cycles $Z \in C_{n-2} (X)$. 
In this case, we have 
\begin{equation*}
\Pic^{1}(X) \xra[\cong]{\overline{\varphi}} 
\J^{n-1}(X) \cong \Hom (\mathcal{H}^{n-1} (X, \Z), \numbersRmodZ)\ ,  
\end{equation*}
and the Abel--Jacobi map 
$\overline{J} \colon \M_{n-2}^{\circ} (X) \cong \J^{n-1}(X)$
is an isomorphism. 
The moduli space $\M_{n-2} (X)$ of Abel $1$--gerbes is generated by 
cycles $Z \in C_{n-2} (X)$ and is given by 
$\M_{n-2} (X) \cong \H_{\D}^{2} (X, \Z)$; that is, 
the moduli space of complex line bundles with unitary connection 
and harmonic curvature. 

\end{example}

\begin{example}
  \label{example3} 
  ~$X = \mathbb{T}^n,  n\geq2, ~d=0,\dots,n-1 \colon$
  
\noindent
  The Jacobi map $\overline{J} \colon \M_{d}^{\circ} (X) \to \J^{d+1} (X)$ 
  is an isomorphism; so is the Picard map 
  $\overline{\alpha} \colon \M_{d}^{\circ} (X) \to \Pic^{n-d-1} (X)$. 
  Here, we take $X = \mathbb{T}^n$, an $n$--dimensional torus with the 
  flat (invariant) Riemannian metric. 
  In this case, the dimension of the Picard--, 
  respectively the Jacobi torus is 
  $\dim H^{d+1} (X, \R) = \binom{n}{d+1}, ~d=0,\dots,n-1$. 
  There is an orthonormal  basis $\{\theta_1, \ldots, \theta_n\}$ of 
  integral, harmonic, invariant $1$--forms which form a framing of 
  the cotangent bundle $T^{\ast} (X)$ and determine an orthonormal 
  basis of $\H^{d+1} (X, \Z)$ by $\theta_{I} = 
  \theta_{i_1} \wedge \dots \wedge \theta_{i_{d+1}}$, where 
  $I = (i_1 < i_2 < \ldots < i_{d+1})$. 
	The dual basis $\{e_1, \ldots , e_n\}$ determines 
	$(d+1)$--subspaces 
	$e_{I} = e_{i_1} \wedge \ldots \wedge e_{i_{d+1}}$ of 
	$\mathbb{R}^n$, respectively basis elements of 
	$\Lambda^{d+1} (\R^n)$ 
	for all multi--indices $I$ as above.  
	By deforming small $(d+1)$--parallelepipeds $P_{I}$ in the 
	direction of $e_{I}$ by the method of Lemma \ref{variation} 
  and considering the families of Jacobi integrals 
  $J_{P_{I}} (\theta_{I}) =  \int_{P_{I}} \theta_{I}$ or 
  their linear combinations, one generates (small) open sets 
  in the range of the Jacobi map. 

\end{example}

\begin{example}
  \label{example4} 
  ~The case $n \geq 2, ~d = 0 \colon$

\noindent
  This is similar to the classical case of 
  divisors on a Riemann surface. In this case, we have 
  \begin{equation*}
\Pic^{n-1}(X) \xra[\cong]{\overline{\varphi}} 
\J^{1}(X) \cong \Hom (\H^{1} (X, \Z), \numbersRmodZ)\ . 
\end{equation*}
  Theorem \ref{modulitheorem} asserts that the 
  Jacobi map $\overline{J} \colon \M_{0}^{\circ} (X) \to \J^{1} (X)$ 
  is an isomorphism. Therefore, so is the Picard map 
  $\overline{\alpha} \colon \M_{0}^{\circ} (X) \to \Pic^{n-1} (X)$. 
  The moduli space $\M_{0} (X)$ of Abel $(n-1)$--gerbes defined by points 
  $\{ p \} \subset X$, whose curvature is the normalized harmonic volume 
  form $\Vol$, satisfies $\M_{0} (X) \cong \H_{\D}^{n} (X, \Z)$;  
  that is, the space of $(n-1)$--gerbes with harmonic curvature. 
  In this case, we have $\M_{0} (X) / \M_{0}^{\circ} (X) \cong H_{0} (X, \Z)$. 
  If $X$ is connected with basepoint $p_0$, the Abel--Jacobi map 
  $\overline{J} \colon B_0 (X) \ra \J^1 (X)$ defines a smooth 
  mapping $j \colon X \ra \J^1 (X)$ by $j (p) ~(\theta) = 
  \int_{p_0}^{p} \theta \mod \Z$, for $\theta \in \H^{1} (X, \Z)$. 
  In turn, the mapping $j$ determines the Abel--Jacobi map $\overline{J}$ 
  completely. To see this, choose regular $1$--simplices $\Gamma_i$, such 
  that $\partial\Gamma_i = \{ p_i \} - \{ q_i \}, ~i=1,\ldots, m, ~m\geq 1$. 
  Then for $\Gamma = \sum_i \Gamma_i$, we have $J_{\partial\Gamma}~(\theta) 
  = \sum_i ~\int_{\Gamma_i} \theta \equiv 
  \sum_i ~(j(p_i)~(\theta) - j(q_i)~(\theta)) \mod \Z$. 
  Thus Abel's Theorem implies that 
  $\sum_i ~j (p_i) = \sum_i ~j (q_i)$, if and only if 
  the $0$--chains $\sum_i ~\{ p_i \}$ and 
  $\sum_i ~\{ q_i \}$ are linearly equivalent. 
  Note that this argument does not prove surjectivity of $\overline{J}$. 
\end{example}

\begin{example}
  \label{example5} 
  ~Riemann surfaces $X$ of genus $g \geq 1, ~n = 2, ~d = 0$: 
   
\noindent
   In this case, we have 
\begin{equation*}
\Pic^{1}(X) \xra[\cong]{\overline{\varphi}} 
\J^{1}(X) \cong \Hom (\mathcal{H}^{1} (X, \Z), \numbersRmodZ)\ , 
\end{equation*}
and 
$\overline{\alpha} \colon \M_{0}^{\circ} (X) \to \Pic^{1}(X)$, 
respectively 
$\overline{J} \colon \M_{0}^{\circ} (X) \to \J^{1}(X)$, 
correspond to the Picard, respectively the Abel--Jacobi 
map of the Riemann surface. 
The Jacobi integral involves path integrals over $1$--chains 
$\Gamma$ of the form $J_{\Gamma} (\theta) = \int_{\Gamma} \theta$.
The Deligne cohomology $H_{\D}^{2} (X, \Z)$ 
is the moduli space of $1$--gerbes; that is \textit{complex line 
bundles with unitary connection}. 

On a Riemann surface the first cohomology group 
$H^{1} (X, \mathcal{M}_{X}^{\ast})$ vanishes 
(cf. \cite{G}, Ch.7, Theorem 12).        
This is a non--trivial consequence of the Riemann--Roch 
theorem and Serre duality. 
>From the exact cohomology sequence 
\begin{equation*}
0 \ra H^{0} (X, \mathcal{O}_{X}^{\ast}) \ra
H^{0} (X, \mathcal{M}_{X}^{\ast}) \xra{D} H^{0}(X, \mathcal{D}_{X}) 
\xrightarrow{\delta_{\ast}} H^{1} (X, \mathcal{O}_{X}^{\ast}) \ra 
H^{1} (X, \mathcal{M}_{X}^{\ast}) 
\end{equation*}
of the divisor sequence 
\begin{equation*}
0 \rightarrow \mathcal{O}_{X}^{\ast} \rightarrow \mathcal{M}_{X}^{\ast} \xrightarrow{D} 
\mathcal{D_{X}} \rightarrow 0\ , 
\end{equation*}
it follows that 
$H^{0}(X, \mathcal{D}_{X}) \xrightarrow{\delta_{\ast}} 
H^{1} (X, \mathcal{O}_{X}^{\ast})$
is surjective and every holomorphic line bundle is the line 
bundle of a divisor. 
In particular, the divisors of degree zero are mapped onto the 
Picard variety 
\begin{equation*}
\Pic (X) = H^{1} (X, \mathcal{O}_{X}) / H^{1} (X, \mathbb{Z}) 
\subseteq H^{1}(X, \mathcal{O}_{X}^{\ast})\ . 
\end{equation*}
$\Pic (X)$ is a complex torus of $\mathrm{dim}_{\mathbb{C}} \Pic (X) = g$, 
the variety of holomorphic line bundles $\mathcal{L}$ 
with ${\Deg} (\mathcal{L}) = c_{1} (\mathcal{L}) = 0$; 
that is, topologically trivial holomorphic line bundles. 
In our context, the Picard torus $\Pic^1 (X)$ is a real torus 
of dimension $2g$ and the above shows that 
$\M_0^{\circ} (X) \cong \Pic^1 (X)$. 
Our proof of 
$\overline{J} \colon \M_{0}^{\circ} (X) \cong \J^{1}(X)$ 
in Theorem \ref{modulitheorem} is much more 
elementary and closer to the direct generation of all 
holomorphic line bundles via divisors in \cite{G}, Ch.7 (c) and 
the Inversion Theorem in \cite{GH}, Ch. 2.2. 
Thus, we choose suitable $1$--simplices $\Gamma_i$ on cycles
$Y_i$ representing a basis $[Y_i] \in H_1 (X, \Z), ~i=1,\ldots, 2g$ 
and deform their endpoints $p_i = \Gamma_i (0,1)$ along the curves $\Gamma_i$ 
by $\Gamma_i (r) (t_0, t_1) = \Gamma_i (t_0 + (1-r)t_1, r t_1), r \in [0,1]$ 
to the fixed initial points $p_{0, i} = \Gamma_i (1, 0)$. 
In this way one generates an open set in the image of the Abel--Jacobi
map by the functionals $J_{\partial\Gamma (r_1, \ldots, r_{2g})} = 
\int_{\Gamma (r_1, \ldots, r_{2g})}$, where 
$\Gamma (r_1, \ldots, r_{2g}) = \sum_{i=1}^{2g} \Gamma_i (r_i)$. 
The same procedure applies also to Example \ref{example4} for the 
generation of $(n-1)$--gerbes defined by points in $X$. 

\end{example}

\begin{example}
  \label{example6}
  ~The case $n = 3, ~d = 0 \colon$ 

\noindent
This is a special case of Example \ref{example4}; 
so, we have 
  \begin{equation*}
\Pic^{2}(X) \xra[\cong]{\overline{\varphi}} 
\J^{1}(X) \cong \Hom (\mathcal{H}^{1} (X, \Z), \numbersRmodZ)\ . 
\end{equation*}
Theorem \ref{modulitheorem} asserts that the 
Abel--Jacobi map $\overline{J} \colon \M_{0}^{\circ} (X) \cong \J^{1}(X)$
is an isomorphism. 
The moduli space $\M_{0} (X)$ of Abel gerbes defined by points 
$\{ p \} \subset X^3$ and whose curvature is the normalized harmonic 
volume form $\Vol$ is given by $\M_{0} (X) \cong \H_{\D}^{3} (X, \Z)$; 
that is, the moduli space of $2$--gerbes with harmonic curvature. 
In this situation, the Abel--Jacobi map was introduced and Abel's 
theorem proved by Hitchin~\cite{H}, Ch. 3.2 and Chatterjee~\cite{Ch} in 
the context of $2$--gerbes. 
This was one of our motivating examples.  

\end{example}

\begin{example}
  \label{example7}
  ~The case $n > 3, ~d = n - 3 \colon$ 

\noindent
Here we look at Abel $2$--gerbes associated to submanifolds 
$M^{n-3} \subset X^n$ of codimension $3$ or more generally to 
cycles $Z \in C_{n-3} (X)$. 
In this case, we have 
\begin{equation*}
\Pic^{2}(X) \xra[\cong]{\overline{\varphi}} 
\J^{n-2}(X) \cong \Hom (\mathcal{H}^{n-2} (X, \Z), \numbersRmodZ)\ , 
\end{equation*}
Theorem \ref{modulitheorem} asserts that the 
Abel--Jacobi map $\overline{J} \colon \M_{n-3}^{\circ} (X) \to \J^{n-2} (X)$ 
is an isomorphism. Therefore, so is the Picard map 
$\overline{\alpha} \colon \M_{n-3}^{\circ} (X) \to \Pic^{2} (X)$. 
The moduli space $\M_{n-3} (X)$ of Abel $2$--gerbes is generated by 
cycles $Z \in C_{n-3} (X)$ and is given by 
$\M_{n-3} (X) \cong \H_{\D}^{3} (X, \Z)$; that is, 
the moduli space of $2$--gerbes with harmonic curvature. 
However, except possibly for the case $n=4, ~d=1$, it is not clear
whether it would be sufficient to consider only codimension three 
submanifolds $M \subset X$; at any rate, in our proof of Theorem 
\ref{modulitheorem} we need to consider cycles $Z \in C_{n-3} (X)$ 
(cf. Bohr--Hanke--Kotschick~\cite{BHK}). 

For codimension $3$ submanifolds $M^{n-3} \subset X^n$, 
the Abel--Jacobi map was also investigated 
and Abel's theorem proved by Hitchin~\cite{H} and 
Chatterjee~\cite{Ch}, Theorem 6.4.2, in the context of $2$--gerbes. 
Moreover, Hitchin in \cite{H}, Theorem 3.2~ 
proves a moduli theorem for families of special Lagrangian 
$3$--tori in a Calabi--Yau $3$--fold via the Abel--Jacobi map. 
Again, these were motivating examples for the present work. 

\end{example}

\begin{theorem}[Moduli Theorem]
  \label{modulitheorem}
  The following statements are equivalent and hold for any compact 
  connected oriented Riemannian manifold $X$ of dimension 
  $n \geq 2,~d = 0,\dots,n-1 \colon$ 

\paritem{$(1)$} 
  The Picard map $\overline{\alpha} \colon \M_d^{\circ} (X) \to \Pic^{n-d-1}(X)$ 
  is an isomorphism. 
  
\paritem{$(2)$} 
  The Abel--Jacobi map $\overline{J} \colon \M_d^{\circ} (X) \to \J^{d+1}(X)$ 
  is an isomorphism.  
  
\paritem{$(3)$} 
  The mapping $\overline{\Lambda} \colon \M_d (X) \to \H_{\D}^{n-d} (X, \Z)$ 
  is an isomorphism. 
 
\paritem{$(4)$}
  Every equivalence class $[\Lambda]$ of $(n-d-1)$--gerbes in the 
  harmonic Deligne cohomology $\H_{\D}^{n-d} (X, \Z)$ 
  can be realized by a unique (up to linear equivalence) 
  Abel gerbe $\Lambda_{Z}$. 
  
\end{theorem}

\begin{proof}
$(1)$ and $(2)$ are equivalent by diagram \eqref{abeljacobi2}. 
To prove that $(3)$ is equivalent to $(1)$, we observe that
\eqref{DeligneES4}, \eqref{cartesian} and \eqref{Pdual} 
determine a commutative diagram 
\begin{equation}
   \label{moduli}
   \vcenter{
     \xymatrix{
       0 \ar[r] & {\M}_d^{\circ} (X) \ar[r] \ar[d]^{\overline{\alpha}} & 
       {\M}_d (X) \ar[r] \ar[d]^{\overline{\Lambda}} & H_d (X) \ar[r] 
       \ar[d]_{\cong}^{PD} &0   \\
       0 \ar[r] & \Pic^{n-d-1} (X) \ar[r]  & {\H}_{\D}^{n-d}(X,\Z) 
       \ar[r]^{d_*} & H^{n-d} (X, \Z) \ar[r] & 0  
     }
   }
\end{equation}
The result follows now from the $5$--lemma. 
$(4)$ is a restatement of $(3)$. 
Thus it suffices to prove $(2)$. 

We choose an (orthonormal) basis 
$\{ \theta_{1}, \dots, \theta_{k} \}, ~k=\dim H^{d+1} (X, \R)$ 
in the integral lattice $\H^{d+1}(X,\numbersZ) \subset \H^{d+1}(X)$ 
of harmonic forms and observe that the Jacobi vector 
\begin{equation}
  \label{jacobivector1}
J_{\partial\Gamma} = 
\begin{pmatrix}
\int_{\Gamma} \theta_1~,  & \dots  &~, \int_{\Gamma} \theta_k	
\end{pmatrix}
\end{equation}
determines the element $\overline{J}_{\partial\Gamma} \in \J^{d+1} (X)$ 
via the expansion $\theta = \sum_{i=1}^{k} a_i \theta_i~, a_i\in \Z$ 
for any $\theta \in \H^{d+1} (X,\Z)$. 
Next, we choose a dual basis $\{\overline{Y}_i\}$ in the integral lattice 
$j_{\ast} H_{d+1} (X, \Z) \subset H_{d+1} (X, \R)$, represented by cycles 
$Y_{i} \in Z_{d+1} (X), i=1,\dots,k$; that is, we have 
\begin{equation}
  \label{dual1}
\int_{Y_i} \theta_j = \delta_{ij}\ . 
\end{equation}
Now, we write $Y_i = \sum_{\el} \Gamma_{i, \el}$, where the $\Gamma_{i, \el}$ 
are regular $(d+1)$--simplices for a triangulation $K$ of $X$. 
Expanding 
\begin{equation}
  \label{jacobimatrix1}
\det 
\begin{pmatrix}
\int_{Y_i} \theta_j 
\end{pmatrix}
_{(i, j=1,\dots,k)} = 1\ , 
\end{equation}
we see that for every $i=1,\dots,k$, there is an $\el_i$ 
such that for $\Gamma_i = \Gamma_{i, \el_i}$ we have 
\begin{equation}
  \label{jacobimatrix2}
\det 
\begin{pmatrix}
\int_{\Gamma_i} \theta_j 
\end{pmatrix}
_{(i, j=1,\dots,k)} > 0\ . 
\end{equation}
We now deform the regular simplices $\Gamma_i$ according to the 
deformation in the proof of Lemma \ref{variation}; so we consider 
$\Gamma_i (r) \colon \simplex^{d+1} \xra{\phi_r} 
\simplex^{d+1} \xra{\Gamma_i} X, ~r \in [0,1]$; 
that is $\Gamma_i ~(r) = \Gamma_i \circ \phi_{r}$. 
First we have from Lemma \ref{variation} $(5)$
\begin{equation*}
  \lim_{r \downarrow 0} \int_{\Gamma_i (r)} \theta = 0, ~i = 1,\ldots, k\ , 
\end{equation*}
for every $\theta \in \H^{d+1} (X, \Z)$, since $\Gamma_i ~(r)$ 
degenerates to a $d$--simplex as $r \downarrow 0$. This is equivalent to 
\begin{equation}
  \label{limit}
\lim_{r \downarrow 0} J_{\partial\Gamma_i (r)} = 0, ~i = 1,\ldots, k\ . 
\end{equation}

>From \eqref{jacobimatrix2} it follows that we can choose 
$\veps > 0$ sufficiently small, such that the smooth mapping  
\begin{equation}
  \label{jacobimatrix3}
r = (r_1,\ldots,r_k) \mapsto 
\begin{pmatrix}
\int_{\Gamma_i (r_i)} \theta_j 
\end{pmatrix}
_{(i, j=1,\ldots,k)} =
\begin{pmatrix}
\int_{\Gamma_1 (r_1)} \theta_1  & \dots  & \int_{\Gamma_1 (r_1)} \theta_k \\	
\vdots & \vdots & \vdots   \\
\int_{\Gamma_k (r_k)} \theta_1  & \dots  & \int_{\Gamma_k (r_k)} \theta_k
\end{pmatrix}
\end{equation}
has positive determinant for $r_i \in (1-\veps, 1]$. 
Moreover, by passing to a subtriangulation of $K$ if necessary, we can achieve 
that our construction takes place in the interior of the fundamental domain 
in the universal cover $\Hom (\mathcal{H}^{d+1} (X, \Z), \numbersR)$ 
of $\J^{d+1} (X)$. 
Therefore the Jacobi vectors (the row vectors in \eqref{jacobimatrix3}) 
\begin{equation}
  \label{jacobivector2}
J_{\partial\Gamma_i (r_i)} = 
\begin{pmatrix}
\int_{\Gamma_i (r_i)} \theta_1~,  & \dots  &~, \int_{\Gamma_i (r_i)} \theta_k	
\end{pmatrix}
\end{equation}
are linearly independent in $\Hom ({\H}^{d+1} (X, \Z), \numbersR)$ 
for $r_i \in (1-\veps,1]$ and $i=1,\dots,k$. 
Setting $\Gamma (r_1, \ldots , r_k) = \sum_{i=1}^{k} \Gamma_i (r_i)$ 
and taking the linear combination of the Jacobi vectors 
\begin{equation}
  \label{lincomb}
J_{\partial\Gamma (r_1, \ldots , r_k)} = 
\sum_{i=1}^{k} J_{\partial\Gamma_i (r_i)} = \sum_{i=1}^{k} 
\begin{pmatrix}
\int_{\Gamma_i (r_i)} \theta_1~,  & \dots  &~, \int_{\Gamma_i (r_i)} \theta_k 
\end{pmatrix}
\end{equation}
gives a 
mapping 
$\Phi \colon D \to \Hom ({\H}^{d+1} (X, \Z), \numbersR) 
\cong \R^{k}$ defined by 
\begin{equation}
\label{phi}
\Phi \colon D \ni r =(r_1, \ldots, r_k) \mapsto 
J_{\partial\Gamma (r_1, \ldots , r_k)} = 
\sum_{i=1}^{k} J_{\partial\Gamma_i (r_i)} ~.   
\end{equation}
Here $D \subset \R^{k}$ is the hypercube (prism) given by 
$r_i \in [0,1], i=1,\ldots,k$. 
$\Phi$ is continuous on D and 
smooth on the interior $B \subset D$, given by 
$r_i \in (0,1)$, $i=1, \ldots, k$. 

The following lemma asserts that the Jacobian $D \Phi$ has positive 
determinant in the neighborhood of an inner point 
$r_{0} = (r_{0, 1}, \ldots, r_{0, k}) \in B$. By the inverse function 
theorem, $\Phi$ is a local diffeomorphism near $r_{0}$ and therefore 
our theorem follows from Proposition \ref{deftwo} $(3)$. 
\end{proof}

\begin{lemma}
  \label{posdet}
  The Jacobian matrix $D \Phi$ is given by
\begin{equation}
\label{detphi}
D \Phi = 
\left( \pdd{\Phi_j}{r_i} \right)
_{(i, j=1,\ldots,k)} = 
\left( \pdd{\int_{\Gamma_i (r_i)} \theta_j} {r_i} \right)
_{(i, j=1,\ldots,k)}\ .  
\end{equation} 
Further, there exists $r_{0} = (r_{0, 1}, \ldots, r_{0, k}) \in B$, 
such that $\det D \Phi (r_0) > 0$. 
\end{lemma}
 
\begin{proof} 
The form of the Jacobian matrix $D \Phi$ follows immediately from 
the fact that each Jacobi vector $J_{\partial\Gamma_i (r_i)}$ in 
\eqref{jacobimatrix3} and in the definition \eqref{lincomb}, 
\eqref{phi} of $\Phi$ depends only on one variable $r_i$.  
In what follows, we will use this fact repeatedly. 
We now inductively use partial differentiation to 
pass from \eqref{jacobimatrix3} to \eqref{detphi}, 
using just the intermediate value theorem of calculus. 
Writing the matrix in \eqref{jacobimatrix3} as a column 
of Jacobi vectors, we know that the determinant is 
positive for $r_i$ in the indicated region, while it goes to 
$0$ for $r_1 \downarrow 0$ by \eqref{limit}. 
Therefore there is a $r_{0, 1} \in (0,1)$, 
such that  
\begin{equation}
  \label{intermedjacobi1} 
\pdd{}{r_1} \vert_{r_{0, 1}} \det
\begin{pmatrix}
J_{\partial\Gamma_1 (r_1)}    \\	
J_{\partial\Gamma_2 (r_2)}    \\ 
\vdots \\
J_{\partial\Gamma_k (r_k)}
\end{pmatrix} = 
\det \begin{pmatrix}
\pdd{J_{\partial\Gamma_1 (r_{0, 1})}}{r_1}    \\	
J_{\partial\Gamma_2 (r_2)}    \\
\vdots \\
J_{\partial\Gamma_k (r_k)} 
\end{pmatrix} > 0\ . 
\end{equation}
Proceeding inductively, we assume that we have 
$r_{0, 1}, \ldots, r_{0, j-1} \in (0,1), 1 < j\leq k$, 
such that the determinant 
\begin{equation}
  \label{intermedjacobi2} 
\det 
\begin{pmatrix}
\pdd{J_{\partial\Gamma_1 (r_{0, 1})} } {r_{1}}    \\	
\vdots \\
\pdd{J_{\partial\Gamma_{j-1} (r_{0, j-1})} }{r_{j-1}}  \\
J_{\partial\Gamma_j (r_j)}  \\
\vdots \\
J_{\partial\Gamma_k (r_k)}
\end{pmatrix} > 0 
\end{equation}
is positive for $r_j, \ldots, r_k \in (1-\veps,1]$. Since 
\eqref{intermedjacobi2} goes to zero as $r_j \downarrow 0$ 
by \eqref{limit}, there is a $r_{0, j} \in (0,1)$, such that 
\begin{equation}
  \label{intermedjacobi3} 
\pdd{}{r_j} \vert_{r_{0, j}} \det
\begin{pmatrix}
\pdd{J_{\partial\Gamma_1 (r_{0, 1}) } } {r_1}    \\	
\vdots \\
\pdd{J_{\partial\Gamma_{j-1} (r_{0, j-1}) } }{r_{j-1}} \\
J_{\partial\Gamma_j (r_j)}  \\
J_{\partial\Gamma_{j+1} (r_{j+1})}  \\
\vdots \\
J_{\partial\Gamma_k (r_k)}
\end{pmatrix} = 
\det \begin{pmatrix}
\pdd{J_{\partial\Gamma_1 (r_{0, 1})} } {r_{1}}    \\	
\vdots \\
\pdd{J_{\partial\Gamma_{j-1} (r_{0, j-1})} }{r_{j-1}}  \\
\pdd{J_{\partial\Gamma_j (r_{0, j})}}{r_{j}}  \\
J_{\partial\Gamma_{j+1} (r_{j+1})}  \\
\vdots \\
J_{\partial\Gamma_k (r_k)}
\end{pmatrix} > 0\ , 
\end{equation}
with $r_{j+1}, \ldots, r_k$ as above. 
This completes the induction. So for $j=k$ we have 
$r_{0} = ( r_{0, 1}, \ldots, r_{0, k} ) \in B$ such that  
\eqref{intermedjacobi3} is positive at $r_0$. But for 
$j=k$, \eqref{intermedjacobi3} is the determinant of the 
Jacobian \eqref{detphi} and the proof is complete.
\end{proof}


\section{Euler and Thom gerbes
\label{secseven}}

In this section we construct the \textit{Euler gerbe} and the 
\textit{Thom gerbe} of an orthogonal bundle, based on the `gerbe 
approach' in \cite{DK2}. 
In Section \ref{seceight} we will investigate the relationship 
between the Euler--, Thom-- and the Abel gerbe.  
We briefly recall the main properties of the construction of 
\textit{characteristic gerbes} from section 5 of \cite{DK2}~, 
which generalizes the classical constructions of secondary 
characteristic classes and `characters' for connections on 
principal $G$-bundles in terms of simplicial forms. 
For the classical constructions we refer to Kamber--Tondeur~\cite{KT},
Chern--Simons~\cite{CnS}, Cheeger--Simons~\cite{CrS} or Dupont--Kamber~\cite{DK1}.

In the following $p\colon P\to X$ is a smooth principal $G$-bundle, $G$
a Lie-group with only finitely many components and $K\subseteq G$ is
the maximal compact subgroup. 
As in Section 5 of \cite{DK2}, 
we fix an invariant homogeneous polynomial $Q\in I^{n+1}(G)$, $n\geq0$, 
such that one of the following 2 cases occur:

\paritem{Case I:} $Q\in\ker(I^{n+1}(G)\to I^{n+1}(K))$.

\paritem{Case II:} $Q\in I^{n+1}_\Z(G)$, that is, there exists an
integral class $u\in H^{2n+2}(BK,\Z)$ representing the Chern-Weil
image of $Q$ in $H^*(BG,\R)\cong H^*(BK,\R)$.

\medbreak
With this notation the \textit{secondary characteristic class}
associated to $Q$ (case I) or $(Q,u)$ (case II) for a connection $A$
on $P\to X$ is a class
\begin{equation}
\label{chargerbe}
  \begin{aligned}{}
    [\Lambda(Q,A)]&\in H^{2n+2}_\D(X)
    \qquad\hbox{in case I,}\\
    [\Lambda(Q,u,A)]&\in H^{2n+2}_\D(X,\Z)
    \qquad\hbox{in case II.}\\
  \end{aligned}
\end{equation}
Note that the characteristic classes in $H^*_\D(X)$ are defined by global 
forms, whereas the classes in $H^*_\D(X,\Z)$ are defined by simplicial forms. 

\begin{enumerate}

	\item
  The classes in \eqref{chargerbe} are natural with respect to 
  bundle maps and compatible coverings.
	
  \item 
  \textit{Curvature formula} :~ 
  \begin{equation}\label{5.8}
    \begin{aligned}
      d\Lambda(Q,A)&=Q(F_A^{n+1})
      \qquad\hbox{in case I}\\
      d\Lambda(Q,u,A)&=\epsilon^*Q(F_A^{n+1})-\gamma
      \qquad\hbox{in case II}
    \end{aligned}
  \end{equation}
  where $\gamma \in \Omega_\Z(|N\U|)$ represents the
  characteristic class $u(P)$ associated with $u$ 
  and $F_A$ is the curvature of $A$. 

  \item
  If $Q(F_A^{n+1})=0$, then 
  \begin{equation}\label{chargerbe1}
    \begin{aligned}{}
      [~\Lambda(Q,A)~]&\in H^{2n+1}(X,\R)
      \qquad\hbox{in case I}\\
      [~\Lambda(Q,u,A)~]&\in H^{2n+1}(X,\R/\Z)
      \qquad\hbox{in case II},
    \end{aligned}
  \end{equation}  
and
  \begin{equation}
    \label{chargerbe2}
    \beta_*[\Lambda(Q,u,A)]=-u(P)
  \end{equation}
  where $\beta_*\colon H^{2n+1}(X,\R/\Z)\to H^{2n+2}(X,\Z)$ is the
  Bockstein homomorphism.
  
  
\end{enumerate}

We shall now use these classes for the case 
$G = SO (2m), ~\Pf \in I_{\Z}^{m} (SO (2m)) $
the Pfaffian polynomial and $u = e \in H^{2m}(BSO (2m), \Z)$ 
the Euler class.
That is, for $Y$ any smooth manifold and $\pi \colon E \to Y$ a 
$2m$-dimensional oriented vector bundle with Riemannian metric 
and metric connection $A$ we obtain a characteristic class 
$[\Lambda (\Pf, e, A)] \in H_{\D}^{2m} (Y, \Z)$ represented by 
a simplicial form 
\begin{equation*}
\Lambda (\Pf, e, A) \in 
\Omega^{2m-1} (\norm{N\covering})
\end{equation*} 
for a suitable covering 
$\covering$ of $Y$, satisfying 
\begin{equation}
  \label{euler}
d \Lambda (\Pf, e, A) = \epsilon^{\ast} \Pf (F_{A}) - e (E)\ . 
\end{equation}  
We shall refer to this form $\Lambda (\Pf, e, A)$ as 
the \textit{Euler gerbe} associated to $E$. 

As is the case for the primary Euler class $e (E)$, there is an 
alternative definition of $\Lambda (\Pf, e, A)$, using the Thom 
space of $E$: 
 
Let $(\mathbb{B}(E), \mathbb{S}(E))$ denote the ball and sphere bundle
of radius $1$. Then up to the choice of a `bump function' in the radial 
direction, the volume form and connection determine a `canonical' 
representative form $U_{E} \in \Omega_{c}^{2m}(\mathbb{B}(E))$, 
with support inside $\mathbb{S} (E)$, for the Thom class in 
$H^{2m}(\mathbb{B}(E), \mathbb{S} (E))$. The restriction of $U_{E}$ to
a neighborhood of the image of the zero section $s \colon Y \to E$ is 
independent of the choice of the bump function. For $F_{A}$ the curvature 
form of $A$, we have 
\begin{equation}
  \label{thomeuler1} 
  s^{\ast} U_{E} = \Pf (F_{A}) \in \Omega^{2m} (Y)\ . 
\end{equation}
For a suitable open covering $\V$ of $\mathbb{B} (E)$, let 
$\beta_{E} \in \Omega_{\Z}^{2m} (\norm{N \V})$ represent the Thom class
of $E$ in $H^{2m} (\mathbb{B}(E), \mathbb{S}(E), \Z)$; that is, 
$\beta_{E}$ vanishes when restricted to $\norm{N \V \cap \mathbb{S}(E)}$. 
Then there exists a simplicial form $\mu_{E} \in \Omega^{2m-1}(\norm{\V})$, 
also with $\mu_{E}$ vanishing in $\norm{N \V \cap \mathbb{S}(E)}$, such that 
\begin{equation}
  \label{thomeuler2}
  d \mu_{E} = \epsilon^{\ast} U_{E} - \beta_{E}\ . 
\end{equation}
Since $H^{2m-1}(\mathbb{B}(E), \mathbb{S}(E)) = 0$, 
the form $\mu_{E}$ is unique modulo 
$\Omega^{2m-1}_{\Z} (\norm{N \V}) + d \Omega^{2m-2} (\mathbb{B}(E))$ 
and hence the Deligne class 
$[\mu_{E}] \in H_{\D}^{2m} (\mathbb{B}(E), \Z)$ is well--defined. 
We shall call $\mu_{E}$ the \textit{Thom gerbe} of $E$. 

\begin{proposition}
  \label{euler1}
  The Thom gerbe determines the characteristic Euler gerbe by 
  the formula: 
  \begin{equation}
  \label{euler2} 
  [s^{\ast} \mu_{E}] = [\Lambda (\Pf, e, A)] \in H_{\D}^{2m} (Y, \Z)\ . 
  \end{equation}
\end{proposition}

\begin{proof}
In the `universal' case (cf. \cite{DK2}, Proposition 5.3), the differential 
of both sides of the equation is $\epsilon^{\ast} \Pf(F_{A}) -e (E)$ by 
\eqref{euler}, \eqref{thomeuler1} and \eqref{thomeuler2}. 
Hence the result follows from the fact that $H^{2m -1} (BSO(2m) , \R) = 0$. 
\end{proof}

We shall now study the Euler and Thom gerbe in particular for 
$E = \nu_{M} = \nu$, where $\nu_{M} \to M$ is the normal bundle 
of a submanifold $M^d \subset X^n$, which we now assume to be of 
even codimension $n-d = 2m$. 
Here $X$ as usual is a compact oriented Riemannian manifold. 
In this case, we identify $\nu$ with a tubular neighborhood $V$ of 
$M \subset X$ and let $V_{0} = \mathbb{B}(\nu) \subset V$. 

Now both $U_{\nu}$ and $\beta_{\nu}$ define (ordinary, respectively 
simplicial) forms on $V$ with support in $\overline{V}_{0}$ and hence 
$\mu_{\nu} \in \Omega^{2m -1} (\norm{N \V \cap \mathbb{B} (\nu)})$ 
extends (non--canonically) to a simplicial form 
$\widetilde{\mu}_{\nu} \in \Omega^{2m -1}(\norm{N \covering})$ 
for a suitable covering $\covering$ of $X$, extending $\V$ on $V$. 
Hence we have 
\begin{equation}
  \label{extthom} 
  [\widetilde{\mu}_{\nu}] \in H_{\D}^{2m}(X, \Z) ~,
\end{equation}  
which we shall call the \textit{extended Thom gerbe}. 
Once $U_{\nu} \in \Omega^{2m}(V)$ is chosen, 
$[\widetilde{\mu}_{\nu}]$ is well-defined, independent of the 
choice of $\covering$ and the choice of $\mu_{\nu}$. 
But it does depend on the `scaling' of $U_{\nu} \in \Omega^{2m} (V)$. 
This of course is not the case for $\mu_{\nu}$ and  
$[s^{\ast} \mu_{\nu}] \in H_{\D}^{2m} (M)$, 
since $U_{\nu}\vert_{V_{0}}$ has a canonical form.


\section{Comparison of the Abel gerbe and 
the Euler gerbe}
\label{seceight}
Continuing with the situation in Section \ref{secseven} of a 
submanifold $M^d \subset X$, of even codimension, we want to 
compare the Euler gerbe with the Abel gerbe associated to $M$. 
Thus let 
$M \subset V_{0} \subset \overline{V}_{0} \subseteq V$ 
be a tubular neighborhood of $M$ and let 
$U_{\nu} \in \Omega_{c}^{n-d} (V_{0})$ be the `canonical' 
Thom class representative. Extend it to $X$ by 0 outside $V_{0}$ 
(also denoted by $U_{\nu}$ ). With this we can define a 
\textit{topologically trivial} gerbe, called the 
\textit{difference gerbe}
\begin{equation*}
[\tau_{M}] \in H_{\D}^{n-d} (X)
\end{equation*}
as follows: 
\par\noindent
>From the beginning of Section \ref{secfour}, recall that 
$F = F_{Z} \vert_{X - \abs{Z}} = \ast~d \ast H_{Z} 
\vert_{X - \abs{Z}}$  is smooth and satisfies 
\eqref{poisson3}; that is $dF = \eta_{Z}\vert_{X - \abs{Z}}$. 
Triangulate $M \subset X$ and choose the covering $\covering$ 
as in Section \ref{secfour}; choose a partition of unity 
$\{ \varphi_i \}_{i=1,\dots,N}$ subordinate to $\covering$ 
and define \textit{smooth} forms in $\Omega^{n-d}(X)$, 
respectively $\Omega^{n-d-1}(X)$: 
\begin{equation*}
\begin{aligned}
\zeta_{0} &= \sum_{i \leq m} \varphi_i ~\eta_{M} \quad , \quad 
\zeta_{1} = \sum_{i > m} \varphi_i ~\eta_{M}\ ,    \\
\zeta_{2} &= \sum_{i > m} d \varphi_i \wedge F = - 
\sum_{i \leq m} d \varphi_i \wedge F\ ,            \\
  G_{1} &= \sum_{i > m} \varphi_i ~F\ .
\end{aligned}
\end{equation*}
Again, we have 
\begin{equation*}
\begin{aligned}
d G_{1} &= \zeta_1 + \zeta_2  \quad , \quad 
\eta_{M} = \zeta_{0} + \zeta_1\ ,     \\ 
\eta_{M} &= (\zeta_{0} - \zeta_2) + d G_1\ , 
\end{aligned}
\end{equation*}
with $\supp (\zeta_{0} - \zeta_2) \subseteq V$. 
Then $d (\zeta_{0} - \zeta_2) = 0$ and hence 
$\zeta_{0} - \zeta_2 = U_{\nu} +d \lambda$, for 
$\lambda \in \Omega^{n-d-1} (X), ~\supp (\lambda) \subseteq V$. 
Then we put $\tau_{M} = \lambda + G_1 \in \Omega^{n-d-1} (X)$, so that 
$d \tau_{M} = d \lambda + d G_1 = 
(\zeta_{0} - \zeta_2 - U_{\nu}) + d G_1 = \eta_{M} - U_{\nu}$; 
that is, $\tau_{M}$ satisfies 
\begin{equation}
  \label{tauequ}
\begin{aligned}
\tau_{M} &= \lambda + G_1 \in \Omega^{n-d-1} (X)\ ,    \\
d \tau_{M} &= \eta_{M} - U_{\nu}\ . 
\end{aligned}
\end{equation}
Hence we get $[\tau_{M}] \in H_{\D}^{n-d}(X)$, and again this 
is well-defined (even independent of the choice of $\covering$ and 
$\{ \varphi_i \}_{i=1,\dots,N}$) once $U_{\nu}$ is chosen. 
Again, since $U_{\nu} \vert_{V_{0}}$ has a canonical form, 
we have 
that $[\tau_{M} \vert_{V_{0}}] \in H_{\D}^{n-d}(V_{0})$ and 
$[\tau_{M} \vert_{M}] \in H_{\D}^{n-d}(M)$ are well--defined. 
Also note that in a neighborhood of $M$ (say $V_{0}$), we have 
$G_1 \vert_{V_{0}} = 0$, so that $\tau_{M} \vert_{V_{0}} = \lambda$.  

\begin{theorem}
  \label{abeleuler1}
 
\paritem{$(1)$} 
In $H_{\D}^{n-d} (X, \numbersZ)$, we have 
\begin{equation}
\label{diffgerbe1}
[\Lambda_{M}] = [\widetilde{\mu}_{\nu}] + \iota_{\ast} ~[\tau_{M}] 
\end{equation}
and all three are well-defined, except that they depend 
on a `scaling' of $U_{\nu}$. 
Further, the characteristic class of $\Lambda_{M}$ 
and $\widetilde{\mu}_{\nu}$ is 
$[\beta_{M}] = [\beta_{\nu}] \in H^{n-d} (X, \Z)$ 
and $\tau_{M}$ is a topologically trivial gerbe 
with curvature $\eta_{M} - U_{\nu}$, where 
$\eta_{M} \in \H^{n-d} (X, \Z)$. 

\paritem{$(2)$} 
In particular, in $H_{\D}^{n-d} (M, \numbersZ)$, we have 
\begin{equation}
  \label{diffgerbe2}
[\Lambda_{M}\vert_{M}] = [\Lambda(\Pf, e, A)] + \iota_{\ast} ~ 
[\tau_{M}\vert_{M}]\ . 
\end{equation}   


\end{theorem}

\begin{proof}
First we observe that the integral simplicial form 
$\beta_{M} = \beta_{Z_{M}}$, 
representing the Poincar\'e dual of $[M] \in H_d (X)$, 
which was constructed at the beginning of the proof of 
Theorem \ref{abelgerbeone} and in Remark 
\ref{submanifold}, can now be chosen to be the integral 
representative $\beta_{\nu}$ of the Thom class for 
the normal bundle $\nu$. 
Since
\begin{equation*}
  \eta_{M}  - \beta_{M} = \eta_{M}  - \beta_{\nu} = 
  (\eta_{M} - U_{\nu}) + (U_{\nu} - \beta_{\nu})\ , 
\end{equation*} 
we have 
\begin{equation*}
  d \Lambda_{M} = d \tau_{M} + d \widetilde{\mu}_{\nu}\ ; 
\end{equation*}
that is
\begin{equation*}
  d (\Lambda_{M} - \tau_{M} - \widetilde{\mu}_{\nu}) = 
  d (\gamma + F_1 - (\lambda + G_1) - \widetilde{\mu}_{\nu}) = 0\ , 
\end{equation*}
where $\gamma, F_1$ are as in the proof of Theorem 
\ref{abelgerbeone}. 
Now, choosing $\covering$ suitable, we can assume 
that $F_1 = \veps^{\ast} (F)~, ~G_1 = F$ on $X - \overline{V}$ 
and therefore $(F_1 - G_1) \vert_{X - \overline{V}} = 0$. 
Since $\gamma, \lambda$ and $\wti{\mu}_{\nu}$ have 
support inside $V$ and since again 
\begin{equation*}
  H^{n-d-1}(\overline{V}, \partial\overline{V}) \cong 
  H^{n-d-1}(\mathbb{B}(\nu), \mathbb{S}(\nu)) = 0\ ,  
\end{equation*}
we get that 
\begin{equation*}
  (\gamma - \lambda) + (F_1 - G_1) - \widetilde{\mu}_{\nu} \in 
  d \Omega^{n-d-2} (\norm{N\covering})\ . 
\end{equation*}
The theorem is proved. 
\end{proof}

\begin{corollary}
  \label{abeleuler2}
Suppose that $[\Lambda_{M}] = 0 \in H_{\D}^{n-d} (X, \numbersZ)$; that is, 
$M \subset X$ is linearly equivalent to zero. 
Then the Euler gerbe $[\Lambda(\Pf, e, A)] \in H_{\D}^{n-d} (M, \numbersZ)$
is topologically trivial, given by the global gerbe 
$[\Lambda(\Pf, e, A)] = - \iota_{\ast} [\tau_{M}\vert_{M}]$.   

\end{corollary}

\begin{proof}
This follows directly from Theorem \ref{abeleuler1} (2). 
\end{proof}



\vfill
\end{document}